\documentclass[a4paper,11pt]{article}

\usepackage{comment}

\usepackage{xifthen}

\newboolean{showDontforget}   
\newboolean{showBacasable}
\newboolean{showComments}

\setboolean{showDontforget}{true}
\setboolean{showBacasable}{false} 
\setboolean{showComments}{true}

\usepackage[%
  hidelinks, 
  bookmarks, 
  linkcolor  = Fblue1,
  citecolor  = Fblue1,
  filecolor  = Fblue1,
  urlcolor   = Fblue1,
  colorlinks = true, 
  pdftitle   = {A Guided Tour of the Equations of Nonlinear Filtering for Diffusion Processes}, 
  pdfauthor  = {Fabien Campillo and Myriam Corso}
]{hyperref}             

\usepackage{multicol}

\usepackage[utf8]{inputenc} 
\usepackage[T1]{fontenc}    
\usepackage[english]{babel}  
\usepackage{float}
\usepackage{amsmath,amsfonts,amssymb}
\usepackage{latexsym}
\usepackage{wrapfig}
\usepackage{wasysym}
\usepackage{lipsum}  
\usepackage{needspace}
  
\usepackage[plain]{algorithm}
\usepackage{algorithmic}

\usepackage[most]{tcolorbox}

\usepackage[
    type={CC},
    modifier={by-nc-sa},
    version={4.0},
]{doclicense}

\usepackage{makecell}
\usepackage{makeidx}
\usepackage[a4paper,left=2.7cm,right=2.7cm,top=2.5cm,bottom=3.5cm]{geometry}
\usepackage{url}
\usepackage{graphicx} 
   \graphicspath{ {pics/} }
\usepackage{fancybox,exscale,rotate,epsfig}
\usepackage{wrapfig}
\usepackage{pifont}
\usepackage{fontawesome5}
\usepackage{letterspace}
\usepackage{supertabular}
\usepackage{rotating}
\usepackage{enumitem} 
	\setlist[enumerate]{leftmargin=30pt, itemsep=0pt}
	\setlist[itemize]{leftmargin=30pt, itemsep=0pt}
\usepackage{tocloft}
\usepackage{setspace}   
\usepackage{framed}

\usepackage{upquote,textcomp}  
\makeatletter
\let \@sverbatim \@verbatim
\def \@verbatim {\@sverbatim \verbatimplus}
{\catcode`'=13 \gdef \verbatimplus{\catcode`'=13 \chardef '=13 }} 
\makeatother

\usepackage{epsfig}
\usepackage{calc}      
\usepackage{ifthen}    

\usepackage{caption}
\usepackage{colortbl}

\usepackage[table,xcdraw,dvipsnames]{xcolor}

\usepackage{titlesec} 
\titleformat{\part}[hang]
   {\huge\bfseries\itshape\flushleft\color{SecColor}}
   {\thepart\ }
   {10pt}
   {\bfseries}
\titleformat{\section}[hang]
   {\Large\bfseries\flushleft\color{SecColor}}
   {\thesection\ }
   {10pt}
   {\bfseries}
\titleformat{\subsection}[hang]
   {\large\bfseries\flushleft\color{SecColor}}
   {\thesubsection\ }
   {0pt}
   {\bfseries}
\titleformat{\subsubsection}[hang]
   {\bfseries\flushleft\color{SecColor}}
   {\thesubsubsection\ }
   {0pt}
   {\bfseries}

\titlespacing*{\section}
  {0pt}    
  {1.5ex}  
  {0.8ex}  

\titleformat{\paragraph}[block]
  {\bfseries\flushleft\color{SecColor}}
  {}
  {0pt}
  {}
\titlespacing*{\paragraph}
  {0pt}
  {1.2ex plus 0.3ex minus 0.2ex}
  {0.5ex}
  
\makeatletter
\renewcommand\subparagraph{%
  \@startsection{subparagraph}{5}{\parindent}%
  {1.ex \@plus 0.5ex \@minus .2ex}
  {-1em}%
  {\normalfont\normalsize\itshape\color{SecColor}}%
}
\makeatother

\definecolor{PA1}{HTML}{E20416}
\definecolor{PA2}{HTML}{F1BE08}
\definecolor{PA3}{HTML}{63D03A}
\definecolor{PA4}{HTML}{C903B5}
\definecolor{PA5}{HTML}{3EBDDC}

\definecolor{Fblue1}{HTML}{3498db}
\definecolor{Fblue2}{HTML}{2980b9}
\definecolor{Fred1}{HTML}{e74c3c}
\definecolor{Fred2}{HTML}{c0392b}
\definecolor{Fgreen1}{HTML}{1abc9c}
\definecolor{Fgreen2}{HTML}{16a085}
\definecolor{Fgrey1}{HTML}{95a5a6}
\definecolor{Fasphalt1}{HTML}{34495e}
\definecolor{Fasphalt2}{HTML}{2c3e50}

\definecolor{lightgray}{rgb}{.9,.9,.9}
\definecolor{darkgray}{rgb}{.4,.4,.4}
\definecolor{mygray}{rgb}{0.5,0.5,0.5}

\definecolor{tab20_blue}{RGB}{31,119,180}
\definecolor{tab20_lightblue}{RGB}{174,199,232}
\definecolor{tab20_orange}{RGB}{255,127,14}
\definecolor{tab20_lightorange}{RGB}{255,187,120}
\definecolor{tab20_green}{RGB}{44,160,44}
\definecolor{tab20_lightgreen}{RGB}{152,223,138}
\definecolor{tab20_red}{RGB}{214,39,40}
\definecolor{tab20_lightred}{RGB}{255,152,150}
\definecolor{tab20_purple}{RGB}{148,103,189}
\definecolor{tab20_lightpurple}{RGB}{197,176,213}
\definecolor{tab20_brown}{RGB}{140,86,75}
\definecolor{tab20_lightbrown}{RGB}{196,156,148}
\definecolor{tab20_pink}{RGB}{227,119,194}
\definecolor{tab20_lightpink}{RGB}{247,182,210}
\definecolor{tab20_gray}{RGB}{127,127,127}
\definecolor{tab20_lightgray}{RGB}{199,199,199}
\definecolor{tab20_olive}{RGB}{188,189,34}
\definecolor{tab20_lightolive}{RGB}{219,219,141}
\definecolor{tab20_cyan}{RGB}{23,190,207}
\definecolor{tab20_lightcyan}{RGB}{158,218,229}

\definecolor{tab10_blue}{RGB}{31,119,180}
\definecolor{tab10_orange}{RGB}{255,127,14}
\definecolor{tab10_green}{RGB}{44,160,44}
\definecolor{tab10_red}{RGB}{214,39,40}
\definecolor{tab10_purple}{RGB}{148,103,189}
\definecolor{tab10_brown}{RGB}{140,86,75}
\definecolor{tab10_pink}{RGB}{227,119,194}
\definecolor{tab10_gray}{RGB}{127,127,127}
\definecolor{tab10_olive}{RGB}{188,189,34}
\definecolor{tab10_cyan}{RGB}{23,190,207}

\definecolor{deep_blue}{RGB}{76,114,176}
\definecolor{deep_orange}{RGB}{221,132,82}
\definecolor{deep_green}{RGB}{85,168,104}
\definecolor{deep_red}{RGB}{196,78,82}
\definecolor{deep_purple}{RGB}{129,114,179}
\definecolor{deep_brown}{RGB}{147,120,96}
\definecolor{deep_pink}{RGB}{218,139,195}
\definecolor{deep_gray}{RGB}{140,140,140}
\definecolor{deep_olive}{RGB}{204,185,116}
\definecolor{deep_cyan}{RGB}{100,181,205}

\usepackage{tikz}
	\definecolor{carmine}{rgb}{0.59, 0.0, 0.09}
	\definecolor{cadmiumred}{rgb}{0.89, 0.0, 0.13}
	\definecolor{cadmiumgreen}{rgb}{0.0, 0.42, 0.24}
	\definecolor{blush}{rgb}{0.87, 0.36, 0.51}
	\definecolor{cadetblue}{rgb}{0.37, 0.62, 0.63}
	\definecolor{crimson}{rgb}{0.86, 0.08, 0.24}
	\definecolor{lightgray}{rgb}{.9,.9,.9}
	\definecolor{darkgray}{rgb}{.4,.4,.4}
	\definecolor{arsenic}{rgb}{0.23, 0.27, 0.29}
	\definecolor{cinnabar}{rgb}{0.89, 0.26, 0.2}
	\colorlet{ColorImmuable}{Orange}
	\colorlet{ColorMutable}{RubineRed}
	\colorlet{ColorCommandePython}{BrickRed}
	\colorlet{ColorCommandePythonIn}{Black!60}
	\colorlet{ColorCommandePythonOut}{Black!30}
	\definecolor{ColorCommandePythonCom}{rgb}{0.5,0.5,0.5}
	\colorlet{ColorLienInterne}{ForestGreen}
	\colorlet{ColorLienExterne}{NavyBlue}
	\colorlet{ColorCite}{cadetblue}
	\colorlet{ColorAnnotation}{carmine!60}
    \colorlet{SecColor}{carmine}
\usetikzlibrary{
   automata,
   arrows.meta,
   backgrounds,
   calc,
   fit,
   hobby,
   math,
   positioning,
   tikzmark,
   shapes,
   shapes.geometric
   }
\tikzset{>=stealth}

\definecolor{fBlue}  {HTML}{4c72b0}
\definecolor{fOrange}{HTML}{dd8452}
\definecolor{fGreen} {HTML}{55a868}
\definecolor{fRed}   {HTML}{c44e52}
\definecolor{fPurple}{HTML}{8172b3}
\definecolor{fBrown} {HTML}{937860}
\definecolor{fPink}  {HTML}{da8bc3}
\definecolor{fTile}  {HTML}{8c8c8c}

\colorlet{f1}{cadmiumred}
\definecolor{f2}{rgb}{0.06667,0.3451,0.62353}
\colorlet{f3}{cadmiumgreen}

\definecolor{mycolor}{rgb}{0.302, 0.416, 0.596} 

\usepackage{needspace}

\definecolor{encadregrey}{gray}{0.93}

\newcommand{\wiki}[1]{\href{#1}{$^\textup{\color{Fblue2}\rm\scriptsize w}$}}

\newcommand{\hrefpep}[1]{\href{#1}{$^\textup{\color{Fblue2}\rm\tiny PEP}$}}

\newcommand{\rmd}   {{{\textrm{\upshape d}}}}

\newcommand{\eqdef}     {\stackrel{\textup{\tiny def}}{=}}

\newcommand{\NN}{\mathcal{N}}

\usepackage{framed}
\definecolor{shadecolor}{gray}{0.75}

\newlength\tindent
\makeatletter
\renewcommand*{\@fnsymbol}[1]{\ensuremath{\ifcase#1\or *\or **\or \ddagger\or
   \mathsection\or \mathparagraph\or \|\or **\or \dagger\dagger
   \or \ddagger\ddagger \else\@ctrerr\fi}}
\makeatother

\usepackage[flushmargin]{footmisc}

\usepackage{listings}

\lstdefinestyle{python}{
  language=Python,
  deletekeywords={...},
  language=Python, 
  alsoletter={:},
  alsoletter={:\n},
  alsoletter={*},
  morekeywords={as,*,:\n,...,.append} 
}

\lstdefinestyle{shell}{
  xleftmargin=0pt, 
  language=sh,
  deletekeywords={for},
  backgroundcolor=\color{white},    
  delim=[il][\bfseries]{BB},
  otherkeywords={port,uninstall,install},
  numbers=none
}

\lstdefinestyle{ipython}{
  xleftmargin=0pt, 
  backgroundcolor=none,    
  language=Python,
  morecomment=[n][\color{Turquoise}]{\ \ \ ...}{:},
  morecomment=[n][\color{Turquoise}]{In\ [}{]\:},
  morecomment=[n][\color{WildStrawberry}]{Out\ [}{]\:},
  morecomment=[n][\color{WildStrawberry}]{Out\[}{]\:},
  backgroundcolor=\color{white},,
  alsoletter={.},
  alsoletter={*},
  alsoletter={:},
  morekeywords={as,*},
  otherkeywords={:\n,.append,.extend,.pop,.count,.index,.replace,.real,.imag,shuffle,%
  .sort,reverse,.split,.join,.keys,.items,copy,.deepcopy,else:,.iterkeys,.find,.upper,%
  .capitalize,.isalpha,.title,.__iter__},
  numbers=none
}

\makeatletter
\def\verbatim@font{\color{Fblue2}\upshape\ttfamily}
\makeatother

\def\Put(#1,#2)#3{\leavevmode\makebox(0,0){\put(#1,#2){#3}}}

{\footnotesize

\setcounter{tocdepth}{2}

}

\newcommand{\brouillon}[1]{\ifthenelse {\boolean{showComments}} {{\footnotesize\color{brouillon} #1}} {}}

\usepackage[textsize=tiny]{todonotes}
\newcommand{\warningr}[2][noinline]{\ifthenelse {\boolean{showComments}} 
{{\todo[#1, color=yellow,size=\tiny]
  {#2}}} {}}
\newcommand{\warningl}[2][noinline]{\ifthenelse {\boolean{showComments}} 
{{\reversemarginpar\todo[#1, color=yellow,size=\tiny]
  {#2}}} {}}
\newcommand{\fabien}[2][noinline]{\ifthenelse {\boolean{showComments}} 
{{\todo[#1, color=yellow,size=\tiny]
  {#2}}} {}}

\newcommand{\fnote}[1]
{{\mbox{}\\\noindent\color{red}\rule{1cm}{2mm}\hfill  #1 \hfill\rule{1cm}{2mm}}\typeout{---------- #1 ------------}}

\newcommand{\Fnote}[1]
{{\mbox{}\\\noindent\color{red}\rule{\textwidth}{2mm}\\  #1 \\ \rule{\textwidth}{2mm}}\typeout{---------- #1 ------------}}

\newcommand{\dontforget}[1]{\ifthenelse {\boolean{showDontforget}} {\fnote{#1}} {}}
\newcommand{\Dontforget}[1]{\ifthenelse {\boolean{showDontforget}} {\Fnote{#1}} {}}
\newcommand{\comments}[1]{\ifthenelse {\boolean{showComments}} {{\color{gray} #1}} {}}

\usepackage{xpatch}
\xpatchbibdriver{online}{urldate}{}{}{}

\usepackage[
  giveninits=true,
  backend=biber,
  style=authoryear,
  isbn=false,
  doi=false,
  url=true,
  eprint=false,
  sorting=nyt,
  maxbibnames=99,
]{biblatex}

\newcommand{\citeannotation}[1]{%
  \iffieldundef{annotation}
    {}
    {\citefield{#1}{annotation}}%
}

\AtEveryBibitem{%
  \iffieldundef{doi}
    {}
    {%
      \clearfield{url}%
      \clearfield{urldate}%
    }%
}

\DeclareFieldFormat{edition}{%
  \ifinteger{#1}
    {\bibstring{edition}~\mkbibordeditio‌​n{#1}}
    {#1\isdot\ edition}}

\renewbibmacro{in:}{%
  \ifentrytype{article}{}{\printtext{\bibstring{in}\intitlepunct}}}	

\AtEveryBibitem{\clearfield{month}} 
\AtEveryCitekey{\clearfield{month}} %
\AtEveryBibitem{\clearfield{series}}
\DeclareFieldFormat{pages}{#1}          
\DeclareFieldFormat[article]{volume}{#1} 
\usepackage{xpatch}

\xpatchbibmacro{journal+issuetitle}{%
  \setunit*{\addspace}%
  \iffieldundef{series}}
  {%
  \setunit*{\addcomma\space}%
  \iffieldundef{series}}{}{}

\newtoggle{bbx:datemissing}

\renewbibmacro*{date}{\toggletrue{bbx:datemissing}%
}

\renewbibmacro*{issue+date}{%
  \toggletrue{bbx:datemissing}%
  \iffieldundef{issue}
    {}
    {\printtext[parens]{%
       \printfield{issue}}}%
  \newunit}

\newbibmacro*{date:print}{%
  \togglefalse{bbx:datemissing}%
  \printdate}

\renewbibmacro*{chapter+pages}{%
  \printfield{chapter}%
  \setunit{\bibpagespunct}%
  \printfield{pages}%
  \newunit
  \usebibmacro{date:print}%
  \newunit}

\renewbibmacro*{note+pages}{%
  \printfield{note}%
  \setunit{\bibpagespunct}%
  \printfield{pages}%
  \newunit
  \usebibmacro{date:print}%
  \newunit}

\renewbibmacro*{addendum+pubstate}{%
  \iftoggle{bbx:datemissing}
    {\usebibmacro{date:print}}
    {}%
  \printfield{addendum}%
  \newunit\newblock
  \printfield{pubstate}}

\renewbibmacro*{url+urldate}{}

\newbibmacro{string+link}[1]{%
  \iffieldundef{doi}
    {%
      \iffieldundef{url}
        {#1}
        {\href{\thefield{url}}{#1}}%
    }
    {\href{https://doi.org/\thefield{doi}}{#1}}%
}

\DeclareFieldFormat{title}{\usebibmacro{string+link}{\mkbibemph{#1}}}
\DeclareFieldFormat[article]{title}{\usebibmacro{string+link}{\mkbibquote{#1}}}
\DeclareFieldFormat[inproceedings]{title}{\usebibmacro{string+link}{\mkbibquote{#1}}}
\DeclareFieldFormat[incollection]{title}{\usebibmacro{string+link}{\mkbibquote{#1}}}

\renewcommand\bibpagespunct{\ifentrytype{article}{\addcolon}{\addcomma}\space}
\renewbibmacro*{volume+number+eid}{%
  \printfield{volume}
  \printfield[parens]{number}%
  \setunit{\addcomma\space}%
  \printfield{eid}}

\DeclareCiteCommand{\citeannotation}
  {\boolfalse{citetracker}%
   \boolfalse{pagetracker}%
   \usebibmacro{prenote}}
  {\printfield{annotation}}
  {\multicitedelim}
  {\usebibmacro{postnote}}

\DeclareCiteCommand{\tabcite}
  {\usebibmacro{prenote}}
  {\usebibmacro{citeindex}%
   \usebibmacro{cite}}
  {\multicitedelim}
  {\usebibmacro{postnote}}
  
\DeclareCiteCommand{\fcitealt}
  {\usebibmacro{cite:init}%
   \usebibmacro{prenote}}
  {\usebibmacro{citeindex}%
   \usebibmacro{cite:comp}}
  {}
  {\usebibmacro{cite:dump}%
   \usebibmacro{postnote}}

\definecolor{sodium}{HTML}{FAAF30}
\definecolor{potassium}{HTML}{AB77B3}
\colorlet{transition}{CornflowerBlue}
\colorlet{open}{ForestGreen}

\newcolumntype{s}{>{\columncolor{sodium!20}}r}
\newcolumntype{p}{>{\columncolor{potassium!20}}r}

\usepackage{mhchem}
\usepackage{chemformula}
\usepackage{nicematrix}

\DeclareFieldFormat{boldimportant}{{\mkbibbold{#1}}{#1}}

\DefineBibliographyExtras{french}{\restorecommand\mkbibnamefamily}

\usepackage{graphicx} 
\usepackage{fancybox,exscale,rotate,epsfig}
\usepackage{wrapfig}
\usepackage{pifont}
\usepackage{letterspace}
\usepackage{supertabular}
\usepackage{rotating}
\usepackage{enumitem} 
	\setlist{leftmargin=0.4cm}
\usepackage{tocloft}
\usepackage{framed}

\usepackage{upquote,textcomp}  
\makeatletter
\let \@sverbatim \@verbatim
\def \@verbatim {\@sverbatim \verbatimplus}
{\catcode`'=13 \gdef \verbatimplus{\catcode`'=13 \chardef '=13 }} 
\makeatother

\usepackage{cprotect}

\usepackage{MnSymbol,wasysym}

\usepackage{mathtools}

\newcommand{\R}{\mathbb{R}}
\renewcommand{\P}{\mathbb{P}}

\newcommand{\crochet}[1]{\langle #1 \rangle}

\newcommand{\E}{\mathbb{E}}
\renewcommand{\P}{\mathbb{P}}
\newcommand{\tE}{{\tilde{\mathbb{E}}}}
\newcommand{\tP}{{\tilde{\mathbb{P}}}}
\newcommand{\tV}{\tilde{V}}

\newcommand{\law}{\textrm{\upshape law}}

\newcommand{\RR}{\mathcal{R}}
\newcommand{\UU}{\mathcal{U}}
\newcommand{\FF}{\mathcal{F}}
\newcommand{\LL}{\mathcal{L}}
\newcommand{\CCb}{\mathcal{C}_{\textrm{\tiny b}}}

\renewcommand{\AA}{\mathcal{A}}
\newcommand{\BB}{\mathcal{B}}
\newcommand{\CC}{\mathcal{C}}
\newcommand{\DD}{\mathcal{D}}
\newcommand{\EE}{\mathcal{E}}
\newcommand{\MM}{\mathcal{M}}

\newcommand{\PP}{\mathcal{P}}
\newcommand{\YY}{\mathcal{Y}}

\renewcommand{\epsilon}{\varepsilon}

\usepackage{listings}

\lstdefinestyle{python}{
  basicstyle=%
    \ttfamily\small
    \lst@ifdisplaystyle\scriptsize\fi,
  language=Python,
  language=Python, 
  alsoletter={:},
  alsoletter={:\n},
  alsoletter={*},
  deletekeywords={...},
  morekeywords={as,*,:\n,...,.append} 
}

\lstdefinestyle{shell}{
  xleftmargin=0pt, 
  language=sh,
  backgroundcolor=\color{white},    
  delim=[il][\bfseries]{BB},
  deletekeywords={for},
  otherkeywords={port,uninstall,install},
  numbers=none
}

\lstdefinestyle{ipython}{
  basicstyle=%
    \ttfamily\upshape\small
    \lst@ifdisplaystyle\scriptsize\fi,
  xleftmargin=0pt, 
  language=Python,
  morecomment=[n][\color{ColorCommandePythonIn}]{\ \ \ ...}{:},
  morecomment=[n][\color{ColorCommandePythonIn}]{In\ [}{]\:},
  morecomment=[n][\color{ColorCommandePythonOut}]{Out\ [}{]\:},
  morecomment=[n][\color{ColorCommandePythonOut}]{Out\[}{]\:},
  backgroundcolor=\color{white},
  alsoletter={.},
  alsoletter={*},
  alsoletter={**},
  alsoletter={:},
  deletekeywords={...},
  morekeywords={as,*,**,:\n,
    __add__,append,add,
    bytes,
    capitalize,copy,count,conjugate,clear,complex,
    deepcopy,difference,difference_update,difference_,discard,
    else,exit,extend,
    False,find,
    getrefcount,getitem,__getitem__,
    imag,index,isupper,islower,isalpha,items,isinstance,__iter__,iterkeys,
      isalnum,isdigit,isspace,intersection,intersection_update,isdisjoint,
      issubset,issuperset,
    join,
    keys,
    __len__,lower,
    pop,
    real,replace,reverse,remove,
    sin,sort,sorted,shuffle,split,swapcase,symmetric_difference,
      symmetric_difference_update,
    title,True,
    upper,union,update,
    values,
    zfill
    },
  otherkeywords={
    },
  numbers=none
}

\lstdefinestyle{neutral}{
  xleftmargin=0pt, 
  backgroundcolor=\color{white},
  numbers=none
  }

\makeatletter
\newcommand\applyCurrentFontsize
{%
  \let\f@sizeS@ved\f@size%
  \let\f@baselineskipS@ved\f@baselineskip%
  \let\basicstyleS@ved\lst@basicstyle%
  \renewcommand\lst@basicstyle%
  {%
      \basicstyleS@ved%
      \fontsize{\f@sizeS@ved}{\f@baselineskipS@ved}%
      \selectfont%
  }%
}
\makeatother

\lstnewenvironment{iplisting}[1][]{%
    \lstset{aboveskip=0.3\baselineskip,
    belowskip=0.3\baselineskip,
    xleftmargin=5pt,
    caption=,
    style=ipython,#1}%
}{}

\usepackage{tikz}
\usetikzlibrary{automata, arrows.meta, positioning}
 
\makeatletter
\def\verbatim@font{\upshape\ttfamily}
\makeatother

\def\Put(#1,#2)#3{\leavevmode\makebox(0,0){\put(#1,#2){#3}}}

\usepackage{tikz}
\usepackage{siunitx}

\usepackage{pgfplots}

\pgfplotsset{compat=1.18}

\pgfplotsset{
  myplot/.style={
    width=0.75\textwidth,
    height=0.45\textwidth,
	xlabel style={font=\huge},
	ylabel style={font=\huge},
	xticklabel style={font=\huge},
	yticklabel style={font=\huge},
	legend style={font=\huge},
    axis lines=middle,
    axis line style={
        black,
        line width=1pt,
        -{Stealth[length=8pt,width=8pt]}
    },
    grid=both,
    major grid style={gray!50},
    minor grid style={gray!25},
    tick style={black},
    ticklabel style={black},
    xlabel style={black},
    ylabel style={black},
    ymin=0,
    ymax=1.2,
    enlargelimits=true
  }
}

\usetikzlibrary{ matrix,      
                 fit,         
                 positioning, 
                 calc, 
                 arrows,shadows,shapes,
                 decorations.markings,
                 arrows.meta,chains,
                 decorations.pathreplacing
}

\tikzset{latentnode/.style={draw, minimum width=5mm, shape=circle, ultra thick, black},
  dagconn/.style={arrows=->, black, thick},
  plate/.style={draw, shape=rectangle, rounded corners=0.5ex, thin, color=gray!50,
    minimum width=3.1cm, text width=3.1cm, align=right, inner sep=5pt, inner ysep=5pt,
    font=\bfseries, text=gray!50,
    append after command={node[text=gray!50, font=\sffamily\scriptsize, 
                               below left= -3pt of \tikzlastnode.south east] {#1}}}
}

\tikzset{
  barbarrow/.style={ 
     >={Straight Barb[left,length=5pt,width=5pt]}
  },
  strike through/.style={
    postaction=decorate,
    decoration={
      markings,
      mark=at position 0.5 with {
        \draw[blue,-] (-2pt,-2pt) -- (2pt, 2pt);
      }
    }
  }
}

\usepackage{chemfig}

\usepackage[american,siunitx]{circuitikz}
\usepackage{amsthm}

\newtheoremstyle{mythmstyle}
  {1ex}
  {1ex}
  {\itshape}
  {0pt}
  {\color{SecColor}\bfseries}
  {.}
  {0.5em}
  {}

\theoremstyle{mythmstyle}
\newtheorem{theorem}{Theorem}[section]

\newtheorem{assumption}[theorem]{Assumption}
\newtheorem{assumptions}[theorem]{Assumptions}

\newtheorem{remark}[theorem]{Remark}

\newcounter{proofstep}

\newcommand{\circnum}[1]{%
\raisebox{0.2ex}{\tikz[baseline=(char.base)]{
\node[
shape=circle,
fill=black!20!white,
text=black,
inner sep=1.2pt,
font=\scriptsize
] (char) {\sf #1};
}}}

\title{\color{SecColor} \bf A Guided Tour of the Equations\\
 \color{SecColor} \bf  of Nonlinear Filtering for Diffusion Processes}

\begin{document}


\author{\color{SecColor} Fabien Campillo\thanks{MathNeuro project -- Inria Branch at University of Montpellier, Montpellier, France.\\ \href{mailto:Fabien.Campillo@inria.fr}{Fabien.Campillo@inria.fr}\\
These  notes were written during the Master’s internship of Myriam Corso (Università di Bologna), carried out in Montpellier under the supervision of Fabien Campillo.\\
The manuscript was prepared with the assistance of ChatGPT, primarily for language editing and stylistic improvements. The authors remain fully responsible for the content.
} \and \color{SecColor}Myriam Corso\footnotemark[1]}
\date{\it\color{SecColor} Inria Branch at University of Montpellier\\
\it\color{SecColor} Montpellier, France \\[1.5em]
\color{SecColor} April 2026}

\maketitle

\newcommand{\x}{\mathbf{x}}
\renewcommand{\o}{\mathbf{o}}
\newcommand{\y}{\mathbf{y}}
\newcommand{\z}{\mathbf{z}}
\renewcommand{\s}{\mathbf{s}}
\renewcommand{\r}{\mathbf{r}}

\newcommand{\sigmav}{\mathbf{v}}
\newcommand{\sigmas}{\mathbf{s}}
\newcommand{\simiid}    {\stackrel{\textrm{\tiny\upshape iid}}{\sim}}

\newcommand{\Cc}{\mathcal{C}_{\textrm{\!\tiny c}}}
\newcommand{\Cb}{\mathcal{C}_{\textrm{\tiny b}}}
\newcommand{\Bb}{\mathcal{B}_{\textrm{\tiny b}}}

\newcommand{\fenumi}{{\upshape({\itshape i}\/)}}
\newcommand{\fenumii}{{\upshape({\itshape ii}\/)}}
\newcommand{\fenumiii}{{\upshape({\itshape iii}\/)}}
\newcommand{\fenumiv}{{\upshape({\itshape iv}\/)}}
\newcommand{\fenumv}{{\upshape({\itshape v}\/)}}
\newcommand{\fenum}[1]{\upshape({\itshape #1}\/)}

\begin{abstract}
These pedagogical notes provide an introduction to nonlinear filtering for diffusion processes. The objective of nonlinear filtering is to estimate an unobserved stochastic state from partial and noisy observations, and thereby characterize the conditional distribution of the state given the observation history.

After introducing the state–observation model, we develop the main tools of the theory, including Markov semigroups, infinitesimal generators,  and change-of-measure methods. This leads naturally to the Kallianpur–Striebel formula, the Zakai equation for the unnormalized filter, and the Kushner–Stratonovich equation for the normalized filter.

The presentation follows, to a large extent, the classical exposition of Bain and Crisan's \emph{Fundamentals of Stochastic Filtering}, while placing particular emphasis on intuition, motivation, and the interpretation of the main concepts and equations. Many technical arguments are revisited, with additional details and
comments intended to facilitate a first reading of the subject.

These notes are not intended to be an exhaustive account of stochastic filtering. Their aim is instead to provide a guided tour of some of its central ideas, methods, and equations.
\end{abstract}

\clearpage
{\small
\setcounter{tocdepth}{3}
\tableofcontents
}

\clearpage
\section{Introduction}

Nonlinear filtering (NLF) is concerned with the estimation of an unobserved
stochastic signal from partial and noisy observations. More precisely, one
considers a \emph{signal process} $(X_t)_{t\ge0}$, also called state process, which cannot be observed
directly, and an \emph{observation process} $(Y_t)_{t\ge0}$, which carries
only indirect information about the signal. The objective is to describe,
at each time $t$, the conditional distribution of $X_t$ given the
observations collected up to time $t$.

In the present notes, we focus on the continuous-time diffusion setting,
where the state evolves according to a stochastic differential equation
and the observations are corrupted by Gaussian white noise:
\begin{subequations}
\label{eq:model}
\begin{align}
  \label{eq:model-x}
  \rmd X_t &= f(X_t)\,\rmd t + g(X_t)\,\rmd W_t\,,
  \quad X_{0} \sim \mu_{0}
\\
  \label{eq:model-y}
  \rmd Y_t &= h(X_t)\,\rmd t + \rmd V_t\,.
\end{align}
\end{subequations}
Although the diffusion framework requires the use of Itô stochastic calculus, which may appear technically demanding at first, it offers a remarkably clear and elegant formulation of the filtering problem. Indeed, many of the central concepts of filtering theory—conditional distributions, innovations processes, change-of-measure techniques, and stochastic partial differential equations—arise naturally in this setting and reveal their underlying connections in a particularly transparent manner.

A central difficulty of nonlinear filtering lies in the nature of its
fundamental object: the conditional law of the state given the observations,
\[
\pi_t(A)
\eqdef
\P\!\left(X_t\in A |Y_s,\;0\le s\le t\right),
\qquad
A\in\BB(\R^n),
\]
where $\BB(\R^n)$ denotes the Borel
$\sigma$-algebra on $\R^n$.

Since the filter is itself a probability measure, the filtering problem
is inherently infinite-dimensional. More precisely, the nonlinear filter
is a stochastic process $(\pi_t)_{t\ge0}$ taking values in the space
$\PP(\R^n)$ of probability measures on
$\R^n$. Consequently, even when the hidden state $X_t$
evolves in a finite-dimensional space, the filter evolves in an
infinite-dimensional state space.

Except in very special situations, most notably the linear Gaussian case,
which leads to the Kalman--Bucy filter, see Section~\ref{sec.finite.dimensional.filters}, the conditional distribution cannot
be characterized by a finite-dimensional system of equations.

The aim of these notes is to present this circle of ideas in a progressive
and intuitive way. Our goal is not to develop the most general or abstract
theory, but rather to explain the main mechanisms underlying nonlinear
filtering in the diffusion setting.

\bigskip

{\it

As the writing of these notes progressed, we gradually realized that we
were following, almost step by step, the path traced by the excellent
monograph of Alan Bain and Dan Crisan \parencite{bain2008a}. By the time we reached the more delicate
questions of the rigorous proofs of existence and uniqueness for the
Zakai equation, prudence suggested that we place our footsteps carefully
in those of the authors~!

The present notes therefore owe much to that work and remain close to it
in spirit. Readers interested in a rigorous and comprehensive treatment
are encouraged to consult the original reference. Our objective here is
instead to offer a guided and intuitive introduction to the main ideas of
nonlinear filtering for diffusion processes.

}

\section{State-space system}

{\it
\color{black!45!SecColor!60!white}

The purpose of this section is to introduce the state-space framework
underlying the nonlinear filtering problem. A state-space model
combines, within a single mathematical formulation, two coupled
stochastic processes: a state process $(X_t)_{t\geq 0}$, representing the
underlying evolution of the system, and an observation process $(Y_t)_{t\geq 0}$,
describing how information about that evolution becomes available to an
observer over time. The model further specifies the mechanism through
which the observations are generated from the hidden state.

A key feature of this framework is its Markovian structure. The hidden
state process $(X_t)_{t\geq 0}$ is Markovian, and the joint process
$(X_t,Y_t)_{t\geq 0}$ is also Markovian. In contrast, the observation process
$(Y_t)_{t\geq 0}$ need not be Markovian when considered on its own. This allows
state-space models to describe a much broader class of observable
dynamics while preserving the tractable Markov structure of the hidden
state.

Another important advantage of the state-space representation is that it
naturally supports Bayesian inference. The Markov property of the hidden
state makes it possible to update the conditional distribution of the
state recursively as new observations become available, leading to the
filtering equations studied throughout these notes.

}

\bigskip

We consider a diffusion state-space model and introduce assumptions ensuring that the system is well posed on a given filtered probability space $\bigl(\Omega,\FF,(\FF_t)_{t\ge0},\P\bigr)$. 
We assume without loss of generality that:
\[
\FF_t = \sigma(W_s,\;0 \le s \le t) \vee\sigma(V_s,\;0 \le s \le t) \vee \sigma(X_0),
\qquad t \ge 0,
\]
completed and made right-continuous\footnote{i.e. $\FF_t = \bigcap_{\epsilon>0}\sigma(W_s,V_s\;0 \le s \le t+\epsilon) \vee \sigma(X_0)\vee \NN$ where $\NN$ is the class of $\P$-null sets where $\AA \vee \BB$ denotes the smallest
$\sigma$-algebra containing both $\AA$ and $\BB$. This is standard in stochastic analysis: completion ensures that probability-theoretic statements are well defined up to null sets, while right-continuity guarantees that the filtration evolves smoothly in time. These modifications do not alter the underlying problem.}, and $\FF=\bigvee_{t\geq 0}\FF_t$. 

In particular, the
state equation admits a unique strong solution, and the observation
process is well defined. As a consequence, the pair process state/observation
 is a continuous adapted process with suitable
integrability properties.

\subsection{State stochastic differential equation}

We start by the state equation \eqref{eq:model-x}. 
The standard hypotheses concerning that stochastic differential equation are:
\begin{assumptions}[Well-posedness of the state equation]
\label{ass:model-wellposed}
The functions $f:\R^n\to\R^n$ and $g:\R^n\to\R^{n\times m}$ are measurable. Moreover:
\begin{enumerate}[leftmargin=40pt,topsep=0pt,itemsep=0pt]
\item 
$f$ and $g$ are globally Lipschitz, i.e. there exists a constant $L>0$ such that
\[
\|f(x)-f(y)\| + \|g(x)-g(y)\| \le L \,\|x-y\|,
\qquad x,y\in\R^n\,.
\]
\item 
The process $(W_t)_{t\ge0}$  is a standard Brownian motion, independent from  the initial condition $X_0$. We also suppose that $X_0$ is a.s. finite.
\end{enumerate}
\end{assumptions}
Note  that the global Lipschitz condition in Hypotheses~\ref{ass:model-wellposed}-\fenumi\ implies that the coefficients $f$ and $g$ have at most linear growth, i.e. $\|f(x)\| + \|g(x)\| \le C \,(1+\|x\|)$.

Under Assumptions~\ref{ass:model-wellposed}, the stochastic differential equation \eqref{eq:model-x}
admits a unique strong solution $(X_t)_{t\ge0}$ (see \cite[Theorem 5.2.1]{oksendal2013a}).

\paragraph{\it Moment assumptions}

Under Assumptions~\ref{ass:model-wellposed}, standard results on
stochastic differential equations imply that for every $p\ge 2$,
\[
\E\|X_0\|^p<\infty
\quad\Longrightarrow\quad
\E\!\left[\sup_{0\le t\le T}\|X_t\|^p\right]<\infty.
\]
A standard proof relies on Hypotheses~\ref{ass:model-wellposed}, Itô's formula, the Burkholder--Davis--Gundy
inequality, and Gronwall's lemma.

\bigskip

In the nonlinear filtering problem, we will need to define conditional
expectations of the form $\E[X_t | \YY_t]$, as well as
quantities related to the estimation error, which involve second-order
moments. It is therefore important to ensure suitable moment bounds for
the state process with this assumption:
\begin{assumption}[Finite second order moments]
\label{ass:model-finite-moments}
We assume that
$\E\|X_0\|^2 < \infty$.
\end{assumption}

Although second moments are sufficient for the formulation of the
filtering problem, the derivation of the Zakai equation requires
stronger integrability estimates. To this end, we strengthen
Assumption~\ref{ass:model-finite-moments} and will assume the existence of
finite third-order moments; see
Assumption~\ref{ass:model-finite-moments.bis}.

\bigskip

Although second moments are sufficient for the formulation of the
filtering problem, the derivation of the Zakai equation requires
stronger integrability estimates. To this end, we strengthen
Assumption~\ref{ass:model-finite-moments} and assume the existence of
finite exponential moments.

As we shall see, the derivation of the Zakai equation relies on the
following chain of logical dependencies:
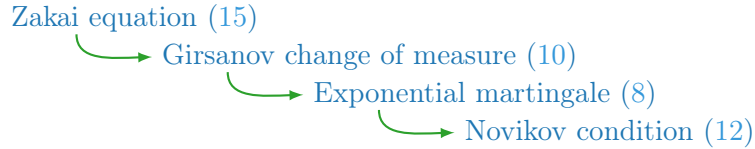
\begin{center}
\begin{tikzpicture}[>=latex]
\node[box,anchor=west,color=tab10_blue]
     (zakai)    at (0,0)
     {Zakai equation \eqref{eq:zakai}};
\node[box,anchor=west,color=tab10_blue]
     (girsanov) at (2,-0.5)
     {Girsanov change of measure \eqref{eq:change.proba}};
\node[box,anchor=west,color=tab10_blue]
     (expo)     at (4,-1)
     {Exponential martingale \eqref{eq:Zt}};
\node[box,anchor=west,color=tab10_blue]
     (novikov)  at (6,-1.5)
     {Novikov condition \eqref{eq.novikov}};
\draw[->,thick,color=tab10_green]
    ($(zakai.west)+(1,-0.2)$)
    to[out=-90,in=180]
    (girsanov.west);
\draw[->,thick,color=tab10_green]
    ($(girsanov.west)+(1,-0.2)$)
    to[out=-90,in=180]
    (expo.west);
\draw[->,thick,color=tab10_green]
    ($(expo.west)+(1,-0.2)$)
    to[out=-90,in=180]
    (novikov.west);
\end{tikzpicture}
\end{center}

When the observation function $h$ is bounded, Novikov's condition
\eqref{eq.novikov} is automatically satisfied. However, this assumption
excludes important situations, including the classical linear--Gaussian
setting. When $h$ is allowed to have linear growth, the verification
of Novikov's condition becomes substantially more delicate and requires
additional integrability properties of the state process.

The exponential moment assumption introduced below is primarily aimed at
providing sufficient conditions for Novikov's criterion. In the bounded case, no such assumption is needed.

\begin{assumption}[Finite exponential moments]
\label{ass:model-exp-moments}
We assume that there exists $\alpha>0$ such that
\begin{equation}
\label{eq.exp.moments.0}
\E\left[
\exp\bigl(\alpha\,\|X_0\|^2\bigr)
\right]
<\infty.
\end{equation}
\end{assumption}

Under this assumption, for every $T>0$, there exists
$\alpha>0$ such that
\begin{equation}
\label{eq.exp.moments.t}
\E\left[
\exp\Bigl(
\alpha\,\sup_{0\le t\le T}\|X_t\|^2
\Bigr)
\right]
<\infty.
\end{equation}

Such estimates can be obtained using exponential versions of Itô's
formula together with standard moment bounds. They hold in particular
for Gaussian initial conditions and provide the exponential
integrability required to verify Novikov-type conditions.

\subsection{Observation equation}

The observation process  \eqref{eq:model-y} can be rewritten as
\begin{align}
  \label{eq:model-y.bis}
  Y_t = \int_0^t h(X_s)\,\rmd s + V_t\,.
\end{align}
This is not a stochastic differential equation, but rather the sum of a
finite variation process and a Wiener process. 

It provides a rigorous
interpretation of the informal model
\begin{align}
  \label{eq:model-y.bis.informal}
  y_t = h(X_t) + \beta_{t}\,,
\end{align}
where $\beta_{t}$ denotes a white Gaussian noise.
In this representation, $h(X_t)$ represents the information about the
hidden state $X_t$ that is accessible through observations. Since the
observation function $h$ is generally not one-to-one, the observations
provide only a partial and potentially ambiguous view of the state,
while $\beta_{t}$ models measurement noise.

The formulation \eqref{eq:model-y.bis} makes \eqref{eq:model-y.bis.informal} precise within the
framework of stochastic calculus. Finally, observing the trajectory
$(Y_s)_{0\le s\le t}$ is equivalent to observing the process
$(y_s)_{0\le s\le t}$ in the informal model. Moreover, the Wiener
process $V_t$ may be viewed as a rigorous mathematical representation
of the integrated white-noise term, since formally
\[
V_t=\int_0^t \dot V_s\,\rmd s.
\]
Although the derivative $\dot V_t$ does not exist as an ordinary
function, this notation provides a useful heuristic interpretation of
$V_t$ as the accumulation of white Gaussian noise.

\begin{assumption}[Well-posedness of the observation equation]
\label{ass:obs-wellposed}
The function $h:\R^n\mapsto\R^d$ is measurable and has at most
linear growth, i.e. there exists $L>0$ such that
\[
\|h(x)\| \le L\,(1+\|x\|), \qquad x\in\R^n\,.
\]

Moreover, $(V_t)_{t\ge0}$ is a $d$-dimensional standard Brownian motion,
independent of $(W_t)_{t\ge0}$ and of the initial condition $X_0$.
\end{assumption}
Under the first part of that hypothesis, and the moment bounds satisfied by $(X_t)_{t\geq 0}$, we have:
\[
\int_0^t \|h(X_s)\|\,\rmd s<\infty\ 
 \text{a.s. for every } t\ge0,
\]
so that $Y_t$ given by \eqref{eq:model-y.bis}
is well defined.

For the derivation of the Zakai equation, this condition is not sufficient
by itself: an additional integrability assumption is required to ensure
that the exponential process \eqref{eq:Zt} arising in the Girsanov transformation is a
true martingale. This will be specified later.

\begin{remark}[Initial observation $Y_0$]
If $Y_0$ is independent of the signal process at time 0, then
conditioning on $Y_0$ does not modify the distribution of the hidden
state. In this case, $Y_0$ contains no useful information for
filtering purposes and may simply be ignored.

If, on the other hand, $Y_0$ depends on the initial state $X_0$,
then $Y_0$ already contains information about the signal. The
filtering problem should therefore start from the conditional
distribution
$\law(X_0|Y_0)$,
rather than from the prior distribution $\law(X_0)$.

In either case, the information carried by $Y_0$ can be absorbed into
the initial condition of the filter. One may therefore assume, without
loss of generality, that
$Y_0=0$.
\end{remark}

\section{Markov semigroup and infinitesimal generator}

{\it
\color{black!45!SecColor!60!white}
The state process $(X_t)_{t\ge0}$ is a Markov diffusion, whose evolution can be described through its associated Markov semigroup and infinitesimal generator. The semigroup characterizes the evolution in time of expectations of functions of the state, while the generator provides a local, differential description of this evolution.

In the context of nonlinear filtering, these objects play a central role. On the one hand, they encode the intrinsic dynamics of the state independently of the observations. On the other hand, they appear explicitly in the evolution equations satisfied by the conditional distribution, such as the Zakai and Kushner–Stratonovich equations. Introducing them allows us to formulate the filtering problem in a framework that naturally connects stochastic processes, conditional expectations, and partial differential equations.
}

\bigskip

The nonlinear filter is a probability measure-valued process. Rather
than studying this object directly, it is often more convenient to
characterize probability measures through their action on suitable test
functions. This naturally leads to a functional viewpoint, in which the evolution of the conditional distribution is characterized through the evolution of expectations of suitable functions of the state.

The starting point of this approach is the Markov diffusion
$(X_t)_{t\ge0}$ itself. Its dynamics are encoded by the associated
Markov semigroup and infinitesimal generator. The semigroup describes
the evolution of expectations of functions of the state, while the
generator provides a local differential description of this evolution.
In this sense, the infinitesimal generator may be viewed as the driving mechanism behind the state process.

Under Assumptions~\ref{ass:model-wellposed}, the process $(X_t)_{t\ge0}$
is a time-homogeneous Markov process (see \cite[Section 7.1]{oksendal2013a}). We define the associated Markov
semigroup $(P_t)_{t\ge0}$ by
\[
  P_t\varphi(x)
  \eqdef
  \E\bigl[\varphi(X_t)\,\big|\, X_0 = x\bigr],
  \qquad x\in\R^n,
\]
for any bounded measurable function $\varphi:\R^n\mapsto\R$.

The family $(P_t)_{t\ge0}$ satisfies the semigroup property:
\[
P_{t+s} = P_t \; P_s\,,\quad t,s\ge0\,,\quad P_0 = \mathrm{Id}\,.
\]
The semigroup describes the evolution of expectations of functions of the
state:
\[
  \E\bigl[\varphi(X_t)\bigr]
  =
  \E\bigl[P_t\varphi(X_0)\bigr].
\]
We now introduce the infinitesimal generator associated with the state
process.

\bigskip

Among the spaces of functions from $\R^n$ to $\R$, we shall
frequently use the following:
\begin{enumerate}
\item
$\Bb(\R^n)$ the space of bounded Borel measurable
functions\footnote{Not to be confused with $\BB(\R^n)$ which denotes the Borel $\sigma$-algebra on $\R^n$.};

\item
$\Cb^2(\R^n)$ the space of twice continuously
differentiable functions whose derivatives up to order two are bounded;

\item
$\Cc^2(\R^n)$ the space of twice
continuously differentiable functions with compact support;
\end{enumerate}
with $\Cc^2(\R^n)\subset \Cb^2(\R^n)\subset \Bb(\R^n)$.

\bigskip

For any test function $\varphi \in \Cb^2(\R^n)$, we define
\begin{tcolorbox}[
 colback=gray!10,colframe=gray!60,width=\textwidth,boxrule=0.3pt,arc=2pt
]
{\vskip-0.4em\hskip-1em\color{gray!80}\sf the infinitesimal generator}
\begin{align}
\label{eq.LL}
  \LL\varphi(x)
  =
  \sum_{i=1}^n f_i(x)\,\partial_i \varphi(x)
  +
  \frac12 \sum_{i,j=1}^n a_{ij}(x)\, \partial_{ij}\varphi(x)
  \,,\qquad x\in\R^n\,,
\end{align}
\end{tcolorbox}
where $a(x):\R^n\to\R^{n\times n}$ is the diffusion matrix defined as
\[
   a(x) \eqdef g(x)\;g(x)^\top\,.
\]
This operator describes the local behavior of the semigroup through the
formal relation
\begin{align}
\label{eq.LL.limit}
  \LL\varphi
  =
  \lim_{t\downarrow 0} \frac{P_t\varphi - \varphi}{t}\,.
\end{align}
This relation suggests the formal notation
\[
\text{`` }P_t=e^{t\LL}\text{ ''}
\]
In other words, one may heuristically think of the semigroup as
\[
\text{`` }P_t
=
\textstyle
\textrm{\upshape Id}
+
t\LL
+
\frac{t^2}{2}\,\LL^2
+
\frac{t^3}{3!}\,\LL^3
+\cdots.
\text{ ''}
\]
This representation is purely formal: for nontrivial diffusion
processes, $\LL$ is typically an unbounded operator, so the
exponential series need not be well defined nor converge in any
classical sense. Nevertheless, the notation $P_t=e^{t\LL}$ provides a
useful heuristic interpretation of the generator as the infinitesimal
counterpart of the Markov semigroup.

\begin{remark}[Domain of the infinitesimal generator]
The domain $\DD(\LL)$ of the infinitesimal generator consists of all
functions for which the limit \eqref{eq.LL.limit} exists. Determining
this domain precisely is often a delicate functional-analytic problem.
Since our objective is not to develop semigroup theory but rather to
derive the equations of nonlinear filtering, we shall avoid these
technical issues and regard the differential operator \eqref{eq.LL} as
the working definition of the generator. Throughout these notes, we
shall therefore work with sufficiently rich families of smooth test
functions, such as $\Cc^2(\R^n)$ or $\Cb^2(\R^n)$, on which the
operator \eqref{eq.LL} is well defined.

In general, the natural domain $\DD(\LL)$ is larger and more intricate
than the smooth function spaces considered here. Since our interest lies
in filtering rather than semigroup theory, we shall not attempt to
characterize it further.
\end{remark}

The infinitesimal generator appears naturally in Itô's formula, indeed:
\[
\varphi(X_t)
=
\varphi(X_0)
+
\int_0^t \LL\varphi(X_s)\,\rmd s
+
\int_0^t \nabla \varphi(X_s)^\top g(X_s)\,\rmd W_s
\]
for any $\varphi \in \Cc^2(\R^n)$.

In particular, the process:
\[
  M_t^\varphi
  \eqdef
  \varphi(X_t)
  - \varphi(X_0) 
  - \int_0^t \LL\varphi(X_s)\,\rmd s
\]
is a martingale.

\bigskip

The generator $\LL$ will play a central role in the filtering
equations. It represents the contribution of the intrinsic dynamics of
the signal process and is the operator governing the Fokker--Planck
equation. In the Zakai equation, this deterministic evolution is
supplemented by terms arising from the observation process.

Working with test functions $\varphi \in \Cc^2(\R^n)$ ensures that
both Itô's formula and the action of the generator $\LL$ are well defined,
which is sufficient for the derivation of the filtering equations.

\section{Filtering problem}

{\it
\color{black!45!SecColor!60!white}
We now introduce the nonlinear filtering problem. Given the observations
of the process $(Y_s)_{0\le s\le t}$, the goal is to estimate the hidden
state $X_t$. This naturally leads to the study of the conditional
distribution of $X_t$ given the observations, which will be described
through conditional expectations with respect to the observation filtration.
}

\bigskip

We first define the observation filtration by
\[
  \YY_t \eqdef \sigma(Y_s,\;0\le s\le t),
  \qquad t\ge0,
\]
completed and made right-continuous.

For every $\varphi\in\Bb(\R^n)$, we define
\begin{tcolorbox}[
 colback=gray!10,colframe=gray!60,width=\textwidth,boxrule=0.3pt,arc=2pt
]
{\vskip-0.4em\hskip-1em\color{gray!80}\sf the nonlinear filter}
\[
  \pi_t(\varphi)
  \eqdef
  \E\left[\varphi(X_t)| \YY_t\right].
\]
\end{tcolorbox}
Since $\varphi$ is bounded, the random variable $\varphi(X_t)$ belongs to
$L^1$, so the conditional expectation is well defined%
\footnote{More precisely, $(\pi_t)_{t\ge0}$ is a version of the regular
conditional distribution of $X_t$ given $\YY_t$, i.e. a
$\YY_t$-measurable random probability measure satisfying
$\pi_t(\varphi)=\E[\varphi(X_t)| \YY_t]$
for all bounded measurable functions $\varphi$. 
In the same way that for two random variables $X_1$ and $X_2$, one can write
$\E[X_1 | X_2] = \psi(X_2)$,
where $\psi(x_2)$ represents the conditional expectation given $X_2 = x_2$,
this function $\psi$ is induced by the regular conditional distribution of
$X_1$ given $X_2$. More precisely, if $\mu(\rmd x_1 | x_2)$ denotes such a version, then
$\psi(x_2)
=
\int x_1\, \mu(\rmd x_1 | x_2)$.}.

\medskip

The process $(\pi_t)_{t\ge0}$ is called the conditional distribution, or
\emph{filter}, of the state given the observations. More precisely, for each
$t\ge0$, $\pi_t$ is a random probability measure on $\R^n$ such that
\[
\pi_t(\varphi)=\int_{\R^n} \varphi(x)\,\pi_t(\rmd x),
\qquad \varphi\in \Bb(\R^n).
\]

A probability measure is uniquely determined by its action on a
measure-determining class of test functions, for instance
$\Bb(\R^n)$ or $\Cb(\R^n)$.  In the derivation of the filtering equations, we shall work with the
smaller class $\Cc^2(\R^n)$, which, as mentioned earlier, is
sufficient for our purposes and is compatible with both Itô's formula
and the infinitesimal generator associated with the state process.

\medskip

\begin{remark}[Smoothing and prediction problems]
Other inference problems, closely related to nonlinear filtering, can
also be considered.

The first is the \emph{smoothing problem}, which consists in computing
\[
[0,T]\ni t \longmapsto \law(X_t | Y_s,\;0\le s\le T).
\]
In contrast to filtering, smoothing exploits observations collected
over the entire time interval $[0,T]$, including those occurring after
time $t$. As a consequence, smoothing is intrinsically an
\emph{offline} problem and cannot be implemented recursively in real
time.

A second problem is the \emph{prediction problem}, which consists in
computing
\[
[T,\infty)\ni t \longmapsto \law(X_t | Y_s,\;0\le s\le T).
\]
Here the objective is to infer the future evolution of the hidden state
beyond the observation horizon $T$.

Both problems are important in applications and can be naturally
derived from the filtering problem. Indeed, once the filter
$\pi_T=\law(X_T| Y_s,\;0\le s\le T)$ is known, it serves as
the starting point for prediction, while smoothing can be obtained by
combining filtering information with information carried by future
observations.
\end{remark}

\needspace{3\baselineskip}
\section{Change of measure and reference probability}

{\it
\color{black!45!SecColor!60!white}
The observation process $Y_t$ contains a drift term $\int_0^t h(X_s)\,\rmd s$ depending on the
unobserved state, which makes the filtering problem intrinsically
nonlinear. A key idea is to remove this dependence by performing a change
of probability measure.

To this end, we introduce a likelihood process $(Z_t)_{t\geq 0}$, constructed as an
exponential martingale, which will serve as the Radon--Nikodym derivative
between the original probability and a new reference probability. Under
this new measure, the observation process becomes a pure noise process,
and the information about the state is transferred to the likelihood
process itself.

This transformation is at the heart of the reference probability method
and will allow us to express the filtering problem in a linear form.
}

\subsection{The likelihood  process $(Z_t)_{t\geq 0}$}

Under the previous assumptions, one readily checks that
\[
\int_0^T \|h(X_s)\|^2\,\rmd s < \infty
\qquad \text{a.s. for every } T>0.
\]
Hence, the process
\[
M_t \eqdef -\int_0^t h(X_s)^\top\,\rmd V_s
\]
is a well-defined continuous local martingale with quadratic variation
\[
\langle M\rangle_t = \int_0^t \|h(X_s)\|^2\,\rmd s
\]
(see, e.g., \cite[p.~35]{oksendal2013a}).
We may therefore define the exponential process
\[
Z_t
\eqdef
\exp\left(
M_t - \frac12 \langle M\rangle_t
\right),
\]
that is,
\[
Z_t
\eqdef
\exp\left(
-\int_0^t h(X_s)^\top\,\rmd V_s
-\frac12\int_0^t \|h(X_s)\|^2\,\rmd s
\right).
\]

\medskip

Using the relation $\rmd V_s = \rmd Y_s - h(X_s)\,\rmd s$, this can be
rewritten as
\begin{tcolorbox}[
 colback=gray!10,colframe=gray!60,width=\textwidth,boxrule=0.3pt,arc=2pt
]
{\vskip-0.4em\hskip-1em\color{gray!80}\sf the likelihood process}
\begin{align}
\label{eq:Zt}
Z_t
\eqdef
\exp\left(
-\int_0^t h(X_s)^\top\,\rmd Y_s
+\frac12\int_0^t \|h(X_s)\|^2\,\rmd s
\right).
\end{align}
\end{tcolorbox}
In the nonlinear filtering context, this latter form is often more
convenient, since it is expressed directly in terms of the observation
process. We also use the notation:
\begin{align}
\label{eq:tidle.Zt}
\tilde Z_t
\eqdef
Z_t^{-1}
=
\exp\left(
\int_0^t h(X_s)^\top\,\rmd Y_s
-\frac12\int_0^t \|h(X_s)\|^2\,\rmd s
\right).
\end{align}
By the general theory of exponential martingales, $(Z_t)_{t\ge0}$ is a
nonnegative local martingale (see \cite[Prop.~3.4, Ch.~IV]{revuz2005a}).

\subsection{Change of probability}

Now we define a new probability measure $\tP$ on $(\Omega,\FF,\P)$ by
\begin{tcolorbox}[
 colback=gray!10,colframe=gray!60,width=\textwidth,boxrule=0.3pt,arc=2pt
]
{\vskip-0.4em\hskip-1em\color{gray!80}\sf the reference probability}
\begin{align}
\label{eq:change.proba}
\left.\frac{\rmd \tP}{\rmd \P}\right|_{\FF_t}
=
Z_t,
\qquad t\ge0\,.
\end{align}
\end{tcolorbox}
The measure $\tP$ is called the \emph{reference probability}
(\footnote{In fact for all $t$ we can define a probability measure $\tP^t$ such that ${\rmd \tP^t}/{\rmd \P}|_{\FF_t}
=
Z_t$, the familly $(\tP^t)_{t\geq 0}$ is consistent, and we can define a probability  $\tP$ equivalent to $\P$ on any $\FF_{t}$, see \cite[p. 55]{bain2008a}.}).

\medskip

The difficulty is that the process $(Z_t)_{t\ge0}$ is in general only a nonnegative local martingale, and hence a supermartingale\footnote{Indeed, if $(\tau_n)_{n\ge1}$ is a localizing sequence, then $Z_{t\wedge\tau_n}$ is a martingale, and the result follows by applying the conditional Fatou lemma as $n\to\infty$.}, but not necessarily a true martingale. In particular, one only has
$\E[Z_t]\le 1$,
so that mass may be lost in the change of probability \eqref{eq:change.proba}. Therefore, it is essential to ensure that $(Z_t)_{t\ge0}$ is a true martingale.

For the moment, the make the following condition:
\begin{assumption}[Well-posedness of the change of measure]
\label{ass:EZ-wellposed}
We assume that
\[
  \E[Z_t]=1,\qquad \text{for all } t\geq 0.
\]
\end{assumption}
This condition is equivalent to requiring that $(Z_t)_{t\ge0}$ is a martingale. We will later provide sufficient conditions ensuring that it holds. Under this assumption, \eqref{eq:change.proba} defines properly the reference probability measure.

\subsection{Girsanov transformation}

With Assumption~\ref{ass:EZ-wellposed}, Girsanov's theorem implies that,
under the reference probability measure $\tP$, the process
$(Y_t)_{t\ge0}$ is a Brownian motion with respect to its natural
filtration $(\YY_t)_{t\ge0}$. Moreover, $(Y_t)_{t\ge0}$ is now independent
of the Brownian motion $(W_t)_{t\ge0}$. 
For a detailed treatment of Girsanov's theorem, see \cite[Chapter VIII]{revuz2005a}.

In fact, under the reference probability measure $\tP$, the state-observation system takes the form
\begin{align*}
  \rmd X_t &= f(X_t)\,\rmd t + g(X_t)\,\rmd W_t,
\\
  \rmd Y_t &= \rmd \tV_t,
\end{align*}
where the state equation is unchanged, and the observation process reduces to a Brownian motion $(\tV_t)_{t\ge0}$, independent of $(W_t)_{t\ge0}$.

\medskip

Under the reference probability $\tilde{\mathbb P}$,
the observation process no longer contains a drift depending on the
state: it becomes a pure noise process that carries no direct
information about the signal. The fact that the ``observation'' process no longer contains any information
about the state is quite paradoxical. In fact, under the reference
probability measure, all information about the signal is transferred to
the likelihood process, which therefore carries the entire dependence
between the state and the observations.

To understand this transformation, it is useful to compare the two
probability measures. Under the original measure $\mathbb P$, the
observation process carries information about the signal through its
drift: the signal drives the observation dynamics, and the model has a
direct physical interpretation. Under the reference probability
$\tilde{\mathbb P}$, by contrast, the observation process is a
Brownian motion independent of the signal. The signal is no longer
visible in the dynamics of $(Y_t)_{t\ge0}$, and its influence is instead transferred
to the likelihood weights.

\begin{figure}
\begin{center}
\begin{tikzpicture}[
  scale=0.8,
  transform shape,
  >=latex,
  box/.style={
    draw,
    thick,
    rounded corners,
    minimum width=6cm,
    minimum height=3.5cm
  },
  every node/.style={font=\footnotesize}
]

\node[box, draw=none, fill=tab10_blue!20]  (Pbox)  at (0,0) {};
\node[box, draw=none, fill=tab10_green!20] (Tbox)  at (8,0) {};

\node[tab10_blue]  at (0,1.25)
{\textbf{original probability \(\P\)}};

\node[tab10_green] at (8,1.25)
{\textbf{reference probability \(\tilde\P\)}};


\node[anchor=west] at (-2.65,0.55)
{\(\rmd X_t=f(X_t)\,\rmd t+g(X_t)\,\rmd W_t\)};

\node[anchor=west] at (-2.65,-0.10)
{
$\tikz[baseline=(dY.base)]{
\node[
rounded corners,
fill=tab10_red!35,
inner sep=2pt
] (dY)
{\(\rmd Y_t=h(X_t)\,\rmd t\)};
}
+\rmd V_t$
};

\node[anchor=west] at (-2.65,-0.95)
{
\(
\pi_t(\varphi)
=
\E\left[
\varphi(X_t)
\mid
\tikz[baseline=(YY.base)]{
\node[
rounded corners,
fill=tab10_red!35,
inner sep=2pt
] (YY)
{\(Y_s,\;s\le t\)};
}
\right]
\)
};


\node[anchor=west] at (5.35,0.55)
{\(\rmd X_t=f(X_t)\,\rmd t+g(X_t)\,\rmd W_t\)};

\node[anchor=west] at (5.35,-0.10)
{\(\rmd Y_t=\rmd\tilde V_t\)};

\node[anchor=west] at (5.35,-0.95)
{
$\rho_t(\varphi)
=
\tilde\E\left[
\varphi(X_t)
\tikz[baseline=(z.base)]{
\node[
rounded corners,
fill=tab10_red!35,
inner sep=2pt
] (z) {$Z_t$};
}
\mid
Y_s,\;s\le t
\right]$
};


\node[
align=center,
color=tab10_red
] (Q) at (4,-3)
{
\it Where is the information about the state?
};


\draw[->, thin, color=tab10_red!35]
(Q.north)
to[out=90,in=-10]
(-.7,-0.20);

\draw[->, thin, color=tab10_red!35]
(Q.north)
to[out=120,in=-90]
(0.95,-0.8);

\draw[-, thin, color=tab10_red!35]
(Q.north)
to[out=40,in=-90]
(8.1,-1.1);

\end{tikzpicture}
\end{center}
\caption{The Girsanov transformation relocates information about the hidden. Under the original probability measure
$\P$, information is carried by the observation process. Under the
reference probability measure $\tilde\P$, the observations does not depend  of the state and the same information about the state is entirely encoded in
the likelihood process $Z_t$.}
\label{fig.reference}
\end{figure}
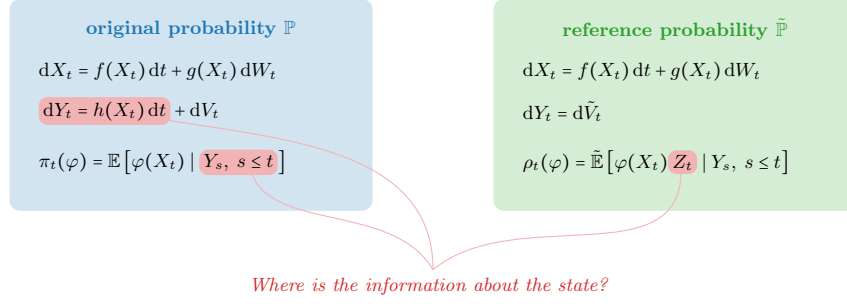

In this way, the information about the signal is shifted from the
dynamics of the observation under $\mathbb P$ to the likelihood process
under $\tilde{\mathbb P}$. \emph{This transformation decouples the
observation noise from the signal and allows one to express the filtering
problem in terms of weighted conditional expectations}, see Figure~\ref{fig.reference}. In particular, it
leads to a linear formulation, which is at the core of the Zakai
equation.

\subsection{Kallianpur--Striebel formula}

This change of measure allows us to express the conditional distribution
of the state in terms of expectations under $\tP$. More precisely,
for any bounded measurable function $\varphi:\R^n\to\R$, we have
\begin{tcolorbox}[
 colback=gray!10,colframe=gray!60,width=\textwidth,boxrule=0.3pt,arc=2pt
]
{\vskip-0.4em\hskip-1em\color{gray!80}\sf the Kallianpur--Striebel formula}
\begin{equation}
\label{eq:kallianpur.striebel}
\pi_t(\varphi)
=
\frac{
\tE\left[ \tilde Z_t\,\varphi(X_t)\,\middle|\, \YY_t \right]
}{
\tE\left[ \tilde Z_t\,\middle|\, \YY_t \right]
}.
\end{equation}
\end{tcolorbox}
The Kallianpur--Striebel formula \eqref{eq:kallianpur.striebel} is essentially a version of
Bayes' rule under a change of measure.

We aim to identify the conditional expectation
$\E[\varphi(X_t)| \YY_t]$. By definition, this is a $\YY_t$-measurable
random variable $\Psi$ such that
\[
\E[H\,\varphi(X_t)] = \E[H\,\Psi]
\]
for every bounded $\YY_t$-measurable random variable $H$.

\medskip

From \eqref{eq:change.proba}, we have
\[
\left.\frac{\rmd \P}{\rmd \tP}\right|_{\FF_t}=\tilde Z_t,
\]
so that
\[
\E[H\,\varphi(X_t)]
=
\tE[H\, \tilde Z_t\varphi(X_t)].
\]
Conditioning under $\tP$ with respect to $\YY_t$, and using the
fact that $H$ is $\YY_t$-measurable, yields
\[
\E[H\,\varphi(X_t)]
=
\tE\bigl[
H\,\tE[\tilde Z_t\varphi(X_t)| \YY_t]
\bigr].
\]

Similarly,
\[
\E[H]
=
\tE[H\, \tilde Z_t]
=
\tE\bigl[
H\,\tE[\tilde Z_t| \YY_t]
\bigr]\,.
\]

\medskip

Comparing these identities, we obtain
\[
\Psi
=
\frac{
\tE[\tilde Z_t\varphi(X_t)| \YY_t]
}{
\tE[\tilde Z_t| \YY_t]
},
\]
which proves \eqref{eq:kallianpur.striebel}.

\paragraph{\it The unnormalized filter}

It is therefore natural to introduce the \emph{unnormalized filter}
$(\rho_t)_{t\ge0}$ defined by
\begin{tcolorbox}[
 colback=gray!10,colframe=gray!60,width=\textwidth,boxrule=0.3pt,arc=2pt
]
{\vskip-0.4em\hskip-1em\color{gray!80}\sf the unnormalized filter}
\[
\rho_t(\varphi)
\eqdef
\tE\left[ \tilde Z_t\,\varphi(X_t)\,\middle|\, \YY_t \right].
\]
\end{tcolorbox}
With this notation, the Kallianpur--Striebel formula simply reads
\[
\pi_t(\varphi)
=
\frac{\rho_t(\varphi)}{\rho_t(1)}.
\]

\subsection{Novikov's condition}
\label{sec.novikov}

Assumption~\ref{ass:EZ-wellposed} is not straightforward to verify, as it
requires sufficient integrability of the exponential of the stochastic
integral defining $(Z_t)_{t\ge0}$. In particular, the likelihood process
is in general only a local martingale, and ensuring that it is a true
martingale requires additional integrability assumptions.

A classical sufficient condition ensuring that
Assumption~\ref{ass:EZ-wellposed} holds is:
\begin{tcolorbox}[
 colback=gray!10,colframe=gray!60,width=\textwidth,boxrule=0.3pt,arc=2pt
]
{\vskip-0.4em\hskip-1em\color{gray!80}\sf the Novikov's condition}
\begin{equation}
\label{eq.novikov}
\E\exp\left(\frac12\int_0^T \|h(X_s)\|^2\,\rmd s\right)<\infty,
\qquad \text{for every } T>0 
\end{equation}
\end{tcolorbox}
(see \cite[Section 3.5.D]{karatzas1991b} or
\cite[Prop.~1.15, Ch.~VIII]{revuz2005a})

There also exists a weaker sufficient condition, known as Kazamaki's
criterion (see \cite[p.~56]{oksendal2013a}), which can be useful in
situations where Novikov's condition fails.

So, when $h$ is bounded, Novikov's condition \eqref{eq.novikov} is
immediately satisfied. However, this assumption excludes important
situations, in particular the linear--Gaussian setting. 

When $h$ has at most linear growth, the verification of
\eqref{eq.novikov} is more delicate and typically relies on exponential
integrability properties of the state process.

Suppose that $h$ has at most linear growth, i.e.\ that
Assumption~\ref{ass:obs-wellposed} holds. We have already seen that,
under this condition and assuming suitable exponential moment bounds on
the initial condition $X_0$ (Assumption~\ref{ass:model-exp-moments}),
the process $(X_t)_{t\ge0}$ admits exponential moments for all $t>0$
(see \eqref{eq.exp.moments.t}).

We now conclude this section by showing that these exponential
integrability properties imply that Novikov's condition
\eqref{eq.novikov} is satisfied.

Indeed, by the linear growth assumption,
\[
\|h(x)\|^2 \le 2\,L^2\,(1+\|x\|^2),
\]
hence
\[
\frac12\int_0^T \|h(X_s)\|^2\,\rmd s
\le
L^2 \,T + L^2 \int_0^T \|X_s\|^2\,\rmd s.
\]
Therefore
\[
\E\exp\left(\frac12\int_0^T \|h(X_s)\|^2\,\rmd s\right)
\le
e^{L^2 \,T}\,
\E\exp\left(L^2 \int_0^T \|X_s\|^2\,\rmd s\right).
\]

Using Jensen's inequality,
\[
\exp\left(L^2 \int_0^T \|X_s\|^2\,\rmd s\right)
=
\exp\left(L^2 \,T \, \frac1T \int_0^T \|X_s\|^2\,\rmd s\right)
\le
\frac1T \int_0^T \exp\bigl(L^2\, T \|X_s\|^2\bigr)\,\rmd s.
\]
Taking expectations yields
\[
\E\exp\left(L^2 \int_0^T \|X_s\|^2\,\rmd s\right)
\le
\frac1T \int_0^T \E\exp\bigl(L^2 \,T \|X_s\|^2\bigr)\,\rmd s,
\]
which is finite by assumption for $\eta \ge L^2 \,T$.

\section{The Zakai equation}

{\it
\color{black!45!SecColor!60!white}

We now derive the Zakai equation, which describes the evolution of the
unnormalized conditional distribution of the signal. Although the
quantity of primary interest is the normalized filter, the unnormalized
formulation plays a central role because it leads to a linear stochastic
evolution equation. This linearity makes the Zakai equation considerably
more tractable than the nonlinear Kushner--Stratonovich equation and is
one of the main reasons why it occupies a central place in modern
filtering theory.

The key idea is to combine the change-of-measure representation provided
by the Kallianpur--Striebel formula with Itô's formula and then take
conditional expectations under the reference probability measure. We
first present the derivation at a formal level, where the essential ideas
are most transparent, and then discuss the additional arguments required
to make the computation rigorous.

}

\bigskip

Recall that the unnormalized filter is defined by
\[
\rho_t(\varphi)
=
\tE\left[ \tilde Z_t\,\varphi(X_t)\,\middle|\, \YY_t \right].
\]
A natural strategy is therefore to apply Itô's formula to the product
$\tilde Z_t\,\varphi(X_t)$ and then take the conditional expectation with
respect to $\YY_t$ under $\tP$.

However, this argument is only formal. Although both $\tilde Z_t$ and
$\varphi(X_t)$ are semimartingales, which allows one to apply Itô's
formula at a purely algebraic level, this is not sufficient to justify
the subsequent probabilistic manipulations. In particular, the process
$(\tilde Z_t)_{t\ge0}$ is in general only a local martingale and may fail
to be integrable. As a consequence, the stochastic integrals arising in
Itô's formula may also be only local martingales, and not true
martingales.

This distinction is crucial: while Itô's formula holds for
semimartingales, the passage to expectations or conditional expectations
requires true martingale properties. Without sufficient integrability,
one cannot guarantee that these stochastic integrals have zero
expectation, nor justify the interchange of limits and expectations. For
this reason, the above computation must be interpreted as formal, and a
regularization procedure is required to make the argument rigorous.

\medskip

Before addressing these technical issues, let us present the formal
derivation of the Zakai equation. This derivation is quite
straightforward and already contains the essential ideas; the technical
difficulties discussed above merely obscure this underlying simplicity.

\subsection{Formal derivation}

Let $\varphi \in \Cc^2(\R^n)$. Under the reference probability
$\tP$, we apply Itô's formula to the product
$\tilde Z_t\varphi(X_t)$.

From \eqref{eq:tidle.Zt}, Itô's formula gives
\[
  \rmd \tilde Z_t
  =
  \tilde Z_t\, h(X_t)^\top \,\rmd Y_t\,.
\]
Also
\[
  \rmd \varphi(X_t)
  =
  \LL\varphi(X_t)\,\rmd t
  +
  \nabla \varphi(X_t)^\top g(X_t)\,\rmd W_t\,,
\]
and finally
\begin{align*}
  \rmd\bigl(\tilde Z_t\,\varphi(X_t)\bigr)
  &=
  \tilde Z_t\,\rmd \varphi(X_t)
  +
  \varphi(X_t)\,\rmd \tilde Z_t
  +
  \rmd \bigl\langle \tilde Z, \varphi(X)\bigr\rangle_t\,.
\end{align*}
Since $(W_t)_{t\geq 0}$ and $(Y_t)_{t\geq 0}$ are independent under $\tP$, the cross-variation term vanishes, and after integrating from $0$ to $t$, we obtain
\begin{align*}
  \tilde Z_t\,\varphi(X_t)
  =
  \varphi(X_0)
  &+
  \int_0^t \tilde Z_{s}\,\LL\varphi(X_s)\,\rmd s
\\
  &
  +
  \int_0^t \tilde Z_{s}\,\nabla \varphi(X_s)^\top g(X_s)\,\rmd W_s
  +
  \int_0^t \tilde Z_{s}\,\varphi(X_s)\, h(X_s)^\top \,\rmd Y_s\,.
\end{align*}
Taking conditional expectation with respect to the observation
filtration $\YY_t$ yields
\begin{align}
\nonumber
  \tE\left[ \tilde Z_t\,\varphi(X_t)\,\middle|\,\YY_t \right]
  =
  \tE[\varphi(X_0)]
  &+
  \int_0^t 
  \tE\left[ \tilde Z_{s}\,\LL\varphi(X_s)\,\middle|\,\YY_s \right]
  \,\rmd s
\\
\label{eq.zakai.formal.E}
  & 
  +
  \int_0^t 
  \tE\left[
  \tilde Z_{s}\,\varphi(X_s)\,h(X_s)^\top 
  \,\middle|\,\YY_s
  \right]
  \,\rmd Y_s\,.
\end{align}

\subsection{Making the derivation rigorous}

A direct application of Itô's formula to the product
\[
\tilde Z_t\,\varphi(X_t)
\]
under the reference probability measure $\tilde{\mathbb P}$ is
formally straightforward. The difficulty arises when one attempts to
take expectations or conditional expectations of the resulting
expression.

The main difficulty stems from the fact that the likelihood process
$\tilde Z_t$ is typically unbounded. Consequently, products such as
$\tilde Z_t\,\varphi(X_t)$ do not automatically satisfy the
integrability properties required for martingale arguments. The
stochastic integrals arising in Itô's formula are therefore, in
general, only local martingales, and additional work is needed to show
that they are true martingales. Without this step, the use of
expectations and conditional expectations cannot be justified.

\medskip

To overcome these difficulties, \cite{bain2008a} introduce the bounded
approximation
\[
\tilde Z_t^\epsilon
\eqdef
\frac{\tilde Z_t}{1+\epsilon \,\tilde Z_t},
\qquad \epsilon>0.
\]
The process $\tilde Z_t^\epsilon$ is positive and bounded by
$1/\epsilon$, and therefore enjoys improved integrability properties.
In particular, the stochastic integrals appearing in Itô's formula for
$\tilde Z_t^\epsilon\,\varphi(X_t)$ are true martingales. This allows
one to apply Itô's formula together with expectations and conditional
expectations in a rigorous way.

\medskip

Rather than working directly with the time-dependent filtration
$(\YY_t)_{t\ge0}$, it is convenient to consider the $\sigma$-algebra
\[
\YY \eqdef \bigvee_{t\ge0} \YY_t,
\]
generated by the entire observation process. Working with this fixed
$\sigma$-algebra simplifies the use of conditional expectations (see
\cite[p.~57]{bain2008a}).

\medskip

We also strengthen Assumption~\ref{ass:model-finite-moments} as follows:
\begin{assumption}[Finite third-order moment]
\label{ass:model-finite-moments.bis}
We assume that
$\E\|X_0\|^3 < \infty$.
\end{assumption}
Under this additional assumption, together with the previous ones, one
can show that
\begin{align}
\label{eq.zakai.rigorous.1}
\int_{0}^t \bigl|\rho_{s}(\|h\|)\bigr|^2\,\rmd s < \infty
\qquad \tP\text{-a.s.}
\end{align}
(see \cite[p.~65]{bain2008a}).
This condition
is required to ensure that stochastic integrals of the form
\[
  \int_0^t \rho_s(\psi)\,\rmd Y_s
\]
are well-defined. Indeed, since $Y$ is a semimartingale with a
Brownian component, the Itô integral is well defined under condition 
\eqref{eq.zakai.rigorous.1}.

Moreover,  condition \eqref{eq.zakai.rigorous.1} is crucial in justifying the interchange of
conditional expectation and stochastic integration, as well as in
passing to the limit in the approximation procedure. It ultimately
follows from the linear growth of $h$ together with suitable moment
bounds on the signal process.

\medskip

We now state the key identities at the level of the bounded
approximation. For any bounded measurable function $\psi$, define
\[
\rho_t^\epsilon(\psi)
\eqdef
\tE\left[
\tilde Z_t^\epsilon\,\psi(X_t)
\,\middle|\,
\YY_t
\right].
\]
Then, using standard properties of conditional expectation and
stochastic integration (see \cite[p.~57--59]{bain2008a}), one can show
that
\begin{enumerate}
\item
$\tE[U |\YY] = \tE[U | \YY_t]$ for any $\FF_t$-measurable random
variable $U$;

\item
$\tE\left[
\int_0^t \tilde Z_s^\epsilon\,\psi(X_s)\,\rmd s
\,\middle|\,
\YY
\right]
=
\int_0^t \rho_s^\epsilon(\psi)\,\rmd s$;

\item
$\tE\left[
\int_0^t \tilde Z_s^\epsilon\,\psi(X_s)\,\rmd Y_s
\,\middle|\,
\YY
\right]
=
\int_0^t \rho_s^\epsilon(\psi)\,\rmd Y_s$;

\item
$\tE\left[
\int_0^t \tilde Z_s^\epsilon\,\psi(X_s)\,\rmd W_s
\,\middle|\,
\YY
\right]
=
0$.
\end{enumerate}
These identities are justified by the boundedness of
$\tilde Z^\epsilon$, which ensures that all stochastic integrals are
true martingales and that conditional expectations can be interchanged
with integration.

\medskip

Applying Itô's formula to $\tilde Z_t^\epsilon\,\varphi(X_t)$ and
taking conditional expectations, one derives a stochastic evolution
equation satisfied by $(\rho_t^\epsilon)_{t\ge0}$. 
Finally, one shows that, as $\epsilon\to0$,
$\rho_t^\epsilon(\varphi)
\to
\rho_t(\varphi)$,
in an appropriate sense (e.g.\ almost surely, in probability, or in
$L^p$). The corresponding identities for the original process
$\tilde Z$ are then obtained by passing to the limit, using dominated
convergence and martingale convergence arguments.

\subsection{Weak form of the Zakai equation}

Equation \eqref{eq.zakai.formal.E} is in fact the
\begin{tcolorbox}[
 colback=gray!10,colframe=gray!60,width=\textwidth,boxrule=0.3pt,arc=2pt
]
{\vskip-0.4em\hskip-1em\color{gray!80}\sf the (weak) Zakai equation}
\begin{align}
\label{eq:zakai}
\rho_t(\varphi)
&=
\rho_0(\varphi)
+
\int_0^t \rho_s(\LL\varphi)\,\rmd s
+
\int_0^t \rho_s(\varphi\,h)^\top \,\rmd Y_s.
\end{align}
for all $t\geq 0$ and $\varphi \in \Cc^2(\R^n)$,
with initial condition $\rho_0$ equal to
the distribution of $X_0$.
\end{tcolorbox}
By construction, the unnormalized filter $(\rho_t)_{t\ge0}$ is therefore a solution of this equation.

It is a linear stochastic evolution equation for a measure-valued process. More precisely, it is a stochastic differential equation taking values in the space of positive measures on $\R^n$.

We have thus established the Zakai equation and, by construction, proved the existence of a solution.

\subparagraph{Fokker-Planck equation}
{\it When $h \equiv 0$, the observation process carries no information about the state, so the conditional distribution reduces to the law of $X_t$. In this case, the unnormalized filter $\rho_t$ is already normalized and coincides with this law, which we denote by $\mu_t$. The Zakai equation then reduces to
\begin{tcolorbox}[
 colback=gray!10,colframe=gray!60,width=\textwidth,boxrule=0.3pt,arc=2pt
]
{\vskip-0.4em\hskip-1em\color{gray!80}\sf the (weak) Fokker--Planck equation}
\begin{align}
\label{eq:fp}
\mu_t(\varphi)
=
\mu_0(\varphi)
+
\int_0^t \mu_s(\LL \varphi)\,\rmd s,
\end{align}
for all $t\geq 0$ and $\varphi \in \Cc^2(\R^n)$,
with initial condition $\rho_0$ equal to
the distribution of $X_0$.
\end{tcolorbox}
This equation, also called the \emph{forward Kolmogorov equation}, governs the time evolution of the distribution of $(X_t)_{t\ge0}$.
}

\subsection{Strong form of the Zakai equation}

We now give a \emph{formal} derivation of the strong form of the Zakai
equation, see \cite[Ch. 7]{bain2008a} for a rigorous derivation.

Assume that, for each $t\ge0$, the unnormalized filter $\rho_t$ admits a
density with respect to Lebesgue measure, say
\[
\rho_t(\rmd x)=r_t(x)\,\rmd x.
\]
Then, for every test function $\varphi$,
\[
\rho_t(\varphi)
=
\int_{\R^n} \varphi(x)\,r_t(x)\,\rmd x
=
\crochet{\rho_t,\varphi}
=
(r_t,\varphi),
\]
where $\crochet{\rho_t,\varphi}$ denotes the duality pairing between a
measure and a test function, and
$(u,v)$ denotes the usual $L^2$ inner product.

With this notation, the weak Zakai equation \eqref{eq:zakai} becomes
\begin{align*}
(r_t,\varphi)
&=
(r_0,\varphi)
+
\int_0^t (r_s,\LL\varphi)\,\rmd s
+
\int_0^t (r_s,\varphi\,h^\top)\,\rmd Y_s.
\end{align*}

We now introduce the formal adjoint $\LL^*$ of the generator $\LL$,
defined by
\[
(u,\LL\varphi)=(\LL^*u,\varphi),
\]
for all sufficiently smooth and decaying functions $u$ and $\varphi$.
A formal integration by parts\footnote{Using $\varphi \in \Cc^2(\R^n)$ (compact support), integration by parts yields\\
\mbox{}$
\qquad\int_{\R^n} \partial_i[ \varphi(x)]\, f_i(x)\, r_t(x)\,\rmd x
=
- \int_{\R^n} \varphi(x)\, \partial_i[f_i(x)\, r_t(x)]\,\rmd x,
$\\
\mbox{}$
\qquad\int_{\R^n} \partial_{ij}[\varphi(x)]\, a_{ij}(x)\, r_t(x)\,\rmd x
=
\int_{\R^n} \varphi(x)\, \partial_{ij}[a_{ij}(x)\, r_t(x)]\,\rmd x
$.} shows that
\[
\LL^*u(x)
=
-\sum_{i=1}^n \partial_i \bigl(f_i(x)\,u(x)\bigr)
+
\frac12 \sum_{i,j=1}^n \partial_{ij} \bigl(a_{ij}(x)\,u(x)\bigr).
\]
Therefore,
\[
(r_s,\LL\varphi)
=
(\LL^*r_s,\varphi),
\]
and the weak Zakai equation can be rewritten as
\begin{align*}
   (r_t,\varphi)
   =
  (r_0,\varphi)
  + \int_0^t (\LL^*r_s,\varphi)\,\rmd s
  + \int_0^t (r_s\,h^\top,\varphi)\,\rmd Y_s\,.
\end{align*}

Since this identity holds for all test functions $\varphi$, we are led,
at least formally, to
\begin{tcolorbox}[
 colback=gray!10,colframe=gray!60,width=\textwidth,boxrule=0.3pt,arc=2pt
]
{\vskip-0.4em\hskip-1em\color{gray!80}\sf the (strong) Zakai equation}
\begin{align}
\label{eq:zakai.strong}
\rmd r_t(x)
&=
\LL^* r_t(x)\,\rmd t
+
r_t(x)\, h(x)^\top \,\rmd Y_t,
\end{align}
for $t>0$ and $x\in\R^n$, 
with initial condition $r_0$ equal to
the density of the law of $X_0$.
\end{tcolorbox}

Equation~\eqref{eq:zakai.strong} is a stochastic partial differential equation (SPDE) satisfied by the density of the unnormalized filter.

\medskip

This derivation is purely formal. Indeed, several points require
justification: the existence of a density $r_t$, the regularity
needed to interpret $\LL^*r_t$, and the validity of the integration by
parts argument. A rigorous derivation typically requires working in
suitable Sobolev or distribution spaces, and interpreting
\eqref{eq:zakai.strong} as a stochastic partial differential equation in
a weak or variational sense. For a comprehensive exposition, we refer to \cite[Chapter 7]{bain2008a}.

\subparagraph{Fokker--Planck equation (strong form)}
{\it 
Under the same assumptions, and formally assuming that $\mu_t$ admits a
density with respect to Lebesgue measure, say
$\mu_t(\rmd x)=q_t(x)\,\rmd x$, the weak formulation leads to 
\begin{tcolorbox}[
 colback=gray!10,colframe=gray!60,width=\textwidth,boxrule=0.3pt,arc=2pt
]
{\vskip-0.4em\hskip-1em\color{gray!80}\sf the (strong) Fokker--Planck equation}
\begin{align}
\label{eq:fp.strong}
  \partial_t q_t(x)
  =
  \LL^* q_t(x),
\end{align}
for all $t\geq 0$ and $x\in\R^n$, with initial condition $q_0$
equal to the density of the law of $X_0$.
\end{tcolorbox}
This equation is the strong form of the forward Kolmogorov equation
and describes the time evolution of the density of $(X_t)_{t\ge0}$.
}

\section{Uniqueness of the solution to the weak Zakai equation}

{\it
\color{black!45!SecColor!60!white}
We now turn to the uniqueness of solutions to the Zakai equation
\eqref{eq:zakai}. Since the equation is linear, uniqueness is naturally
studied by considering the difference of two solutions and proving that
any solution with zero initial condition must vanish identically.

The proof is based on a duality argument. The idea is to test the
forward measure-valued equation against solutions of a suitably chosen
backward Kolmogorov equation. This choice makes the drift terms cancel
and reduces the problem to a martingale identity. To extract pointwise
information from this identity, one uses a separating family of
exponential martingales adapted to the observation filtration.

We first present the argument at a formal level, in order to highlight
the main mechanism of the proof. We then introduce the functional spaces
and regularity assumptions needed to turn this formal computation into a
rigorous uniqueness theorem.
}

\subsection{A formal uniqueness argument}
\label{sec.zakai.uniqueness.formal}

$\color{SecColor}\blacktriangleright$
The Zakai equation is linear. Therefore, to prove uniqueness, it is
sufficient to show that the only solution with zero initial condition is
the trivial solution.

Indeed, let $\rho^1$ and $\rho^2$ be two solutions of the Zakai
equation with the same initial condition, and define
\[
\nu_t \eqdef \rho_t^1-\rho_t^2.
\]
By linearity, $\nu_t$ satisfies the same equation, but with zero initial
condition:
\begin{equation}
\label{eq.zakai.zakai.ini0}
\rmd \nu_t(\varphi)
=
\nu_t(\LL\varphi)\,\rmd t
+
\nu_t(\varphi\,h^\top)\,\rmd Y_t,
\qquad
\nu_0=0,
\end{equation}
for every admissible test function $\varphi$.

The uniqueness problem is therefore reduced to proving that any solution
of the above equation with zero initial condition must vanish
identically. More precisely, we shall show that
\begin{equation}
\label{eq.zakai.nu.0}
\nu_t(\varphi)=0,
\qquad
\text{for all } t\ge 0
\text{ and all admissible test functions } \varphi.
\end{equation}
Since a measure is uniquely determined by its action on test functions,
it follows that
$\nu_t=0$,
for all $t\ge 0$.
Consequently,
$\rho_t^1=\rho_t^2$,
for all $t\ge 0$
and uniqueness follows.

\bigskip

$\color{SecColor}\blacktriangleright$
To establish \eqref{eq.zakai.nu.0} at any given time $T>0$, we note that for every test function $\varphi$, the quantity
$\nu_T(\varphi)$ is a $\YY_T$-measurable random variable. To prove that
$\nu_T(\varphi)=0$, it is natural to introduce a separating family
$\RR$ of $\YY_T$-measurable bounded random variables. By this we mean that,
for every integrable $\YY_T$-measurable random variable $\xi$,
if 
$\tilde{\E}\left[R\,\xi\right]=0$ for all   $R\in\RR$,
then 
$\xi=0$ $\tilde{\P}$-a.s..
Applying this criterion to
\[
\xi=\nu_T(\varphi),
\]
it is enough to prove that
\[
\tilde{\E}\left[R\,\nu_T(\varphi)\right]=0,
\qquad \forall R\in\RR,
\]
in order to conclude \eqref{eq.zakai.nu.0} for any time $T>0$.

\bigskip

$\color{SecColor}\blacktriangleright$
The idea is then to use the following :
\begin{tcolorbox}[colback=gray!10,colframe=gray!60,width=\textwidth,boxrule=0.3pt,arc=2pt]
{\vskip-0.4em\hskip-1em\color{gray!80}\sf complex exponential
$\YY_t$-martingales separating family}
\begin{equation}
\label{eq.zakai.uniqneness.R_t}
R_t
=
\exp\left(
i\int_0^t r_s^\top\,\rmd Y_s
+
\frac12\int_0^t \|r_s\|^2\,\rmd s
\right)
\end{equation}
where $r$ ranges over the class of bounded deterministic functions
$r:[0,T]\to\R^d$.
\end{tcolorbox}
The resulting family of martingales plays the role
of Fourier modes on the observation path space and forms the separating
family needed in the uniqueness proof. We defer the discussion of the exponential martingale $R_t$ and its
separating property until the end of this section.

The duality trick that follows requires test functions that evolve in
time. We therefore extend the weak form of the Zakai equation from
time-independent test functions $\varphi(x)$ to sufficiently regular
time-dependent test functions
$\phi(t,x)$.

If $(t,x)\to\phi(t,x)$ is sufficiently smooth, then an application of the chain rule yields to the
\begin{tcolorbox}[colback=gray!10,colframe=gray!60,width=\textwidth,boxrule=0.3pt,arc=2pt]
{\vskip-0.4em\hskip-1em\color{gray!80}\sf weak Zakai equation with time-dependent test functions}
\begin{equation}
\label{eq:zakai.t}
  \rmd \nu_t(\phi_t)
  =
  \nu_t\left(
    \partial_t \phi_t+\LL\phi_t
  \right)\rmd t
  +
  \nu_t(\phi_t \,h^\top)\,\rmd Y_t\,,\quad \nu_t(\phi_0)=0\,.
\end{equation}
\end{tcolorbox}
Here and in the sequel, we use the shorthand notation
$\phi_t(x)\eqdef \phi(t,x)$.

Next, observe that a direct application of It\^o's formula shows that
\[
  \rmd R_t = i\,R_t\,r_t^\top\,\rmd Y_t\,.
\]
Applying the product rule gives
\begin{align*}
  &
  \rmd\bigl(R_t\,\nu_t(\phi_t)\bigr)
  =
  R_t\,\rmd\nu_t(\phi_t)
  +
  \nu_t(\phi_t)\,\rmd R_t
  +
  \rmd\langle R,\nu(\phi)\rangle_t\,.
\end{align*}
The quadratic variation term is
\[
  \rmd\langle R,\nu(\phi)\rangle_t
  =
  i\,R_t\,r_t^\top\,\nu_t(\phi_t\, h)\,\rmd t.
\]
because both martingale parts are driven by $\rmd Y_t$.
Hence:
\begin{align}
\label{eq.zakai.uniqneness.R_t.nu_t}
  \rmd\bigl(R_t\,\nu_t(\phi_t)\bigr)
  =
  R_t\;
  \nu_t\left(
    {\color{tab10_green}
    \partial_t\phi_t+\LL\phi_t
    +
    i\,h^\top r_t\,\phi_t
    }
  \right)\rmd t
  +
  R_t\;
  \bigl[
    \nu_t(\phi_t \,h)
    +
    i\,\nu_t(\phi_t)\,r_t
  \bigr]^\top\,\rmd Y_t .
\end{align}

\bigskip

$\color{SecColor}\blacktriangleright$
Since probabilists are rarely happier than when a martingale is in sight,
usually for good reasons, we now bring in the very useful
{\color{tab10_green}duality trick}\footnote{
The duality trick is much older than filtering theory. At its core, it
is the classical PDE idea of pairing a forward equation with a solution
of the corresponding adjoint backward equation in order to cancel drift
terms. What is specific to the Zakai equation is the combination of
this duality argument with exponential martingales, which play the role
of Fourier modes on Wiener space. The resulting method goes back to the
work of Kunita and is often referred to in the filtering literature as
the \emph{Bensoussan method} (see \cite{bensoussan1992b}).
}.

The key idea is to exploit the duality between the forward Zakai
equation and a suitable backward partial differential equation. More
generally, the duality argument consists in testing a forward
measure-valued evolution against a family of functions solving an
appropriately chosen backward equation. The backward dynamics are
designed so that the drift terms cancel in the evolution of the duality
pairing. One is then left with a much simpler identity, typically of
martingale type, shared by all solutions. Since a sufficiently rich
family of test functions uniquely determines a measure, uniqueness
follows.

In the present setting, the trick consists in choosing a suitable class
of time-dependent test functions $\phi_t$ that exactly cancel the
drift term appearing in the dynamics of
$R_t\,\nu_t(\phi_t)$. This is achieved by defining
$\phi_t$ as the solution of the backward equation
\begin{tcolorbox}[colback=gray!10,colframe=gray!60,width=\textwidth,boxrule=0.3pt,arc=2pt]
{\vskip-0.4em\hskip-1em\color{gray!80}\sf backward dual equation}
\begin{equation}
\label{eq.zakai.uniqneness.backward.pde}
\color{tab10_green}
\partial_t \phi_t(x)
+
\LL \phi_t(x)
+
i\,h^\top(x)\,r_t\,\phi_t(x)
=
0,
\qquad
0\le t\le T,
\qquad
\phi_T(x)\equiv\varphi(x)\,,
\end{equation}
where the terminal condition $\varphi$ is completely arbitrary.
\end{tcolorbox}
We
deliberately leave it free, since later in the proof we shall vary
$\varphi$ over a sufficiently rich class of test functions. The
resulting family of backward solutions
$(\phi_t)_{0\le t\le T}$ will allow us to probe the measure
$\nu_T$ in all possible directions and ultimately conclude that
$\nu_T=0$.

The appearance of a backward equation is natural. Our goal is to obtain
information about the quantity $\nu_T(\varphi)$ at the terminal time
$T$. We therefore start from the prescribed terminal condition
$\varphi$ and construct a family of test functions satisfying
$\phi_T=\varphi$. Since the value is specified at the final time
rather than at the initial time, the corresponding PDE must be solved
backward in time. More deeply, this equation is the adjoint counterpart
of the forward evolution satisfied by the measure-valued process
$\nu_t$, which is precisely why it is able to cancel the drift term
in the duality argument.

\bigskip

$\color{SecColor}\blacktriangleright$
Finally as  $R_0\,\nu_0(\phi_0)=0$, we get:
\begin{align*}
  R_t\,\nu_t(\phi_t)
  =
  \int_{0}^t R_s
  \left[ 
    \nu_s(\phi_s\, h^\top)
    +
    i\,\nu_s(\phi_s)\,r_s^\top
  \right]\rmd Y_s \,.
\end{align*}
which is  a local $\YY_{t}$-martingale, and under integrability assumptions a true martingale. Taking expectations yields
\begin{equation}
\label{eq.zakai.uniqneness.separation.cond}
\tilde{\E}
\left[
R_T\,\nu_T(\varphi)
\right]
=
0.
\end{equation}
This identity holds for every terminal test function $\varphi$ and every
choice of the deterministic function $r$. The family of random variables  $R_T$ is sufficiently rich to separate
$\YY_T$-measurable random variables, see later. It follows that
$\nu_T(\varphi)=0$
for every test function $\varphi$, and therefore
$\nu_T=0$.
Since $T$ is arbitrary, we conclude that
\[
\rho_t^1=\rho_t^2,
\qquad \textrm{for all }t\ge0.
\]
Thus, the Zakai equation admits at most one solution.

\bigskip

{\it 
The above argument is purely formal. In particular, we have implicitly
assumed that the solution lives in some natural function space -- so
natural, in fact, that we did not even bother to name it  -- and that
uniqueness holds in that space. One may suspect that this is somewhat
naive!

The first and most important question is therefore to determine the
appropriate  space in which the Zakai equation should be posed and
in which uniqueness can reasonably be established. Only once this issue
has been clarified does the uniqueness statement acquire a precise
 meaning.

These questions will be addressed in the following sections.
}

\subsection{Uniqueness: in which space and under which assumptions on the coefficients?}
\label{sec.zakai.uniqueness.spaces}

{\it

Analysts like to say that, when investigating existence and uniqueness,
once the right function space has been found, the equation is already
half solved. }

\medskip

The question of uniqueness of the Zakai equation cannot be addressed
without first specifying the class of admissible solutions. Indeed,
uniqueness is always relative to a given functional space: one seeks to
determine whether the equation admits a unique solution within a
prescribed class of stochastic processes. This choice is therefore
crucial.

On the one hand, the space must be sufficiently large to contain the
solution arising from the filtering problem. On the other hand, if the
space is chosen too large, uniqueness may fail, or become considerably
more difficult to establish. In particular, enlarging the class of
admissible solutions may allow pathological processes that satisfy the
equation in a weak sense but do not correspond to the filtering problem.

\medskip

The Zakai equation \eqref{eq:zakai} 
describes the evolution of positive measure-valued
processes. To prove uniqueness, one usually considers two candidate
solutions $\rho^1$ and $\rho^2$ and studies their difference
\[
\nu_t \eqdef \rho_t^1-\rho_t^2.
\]
Because the Zakai equation is linear, $\nu$ is itself a solution of
the Zakai equation with zero initial condition.

There is, however, a small catch. While $\rho^1$ and $\rho^2$ are
positive measures, their difference is generally only a signed measure.
The uniqueness problem is therefore transformed into the question of
whether the Zakai equation admits a nontrivial signed measure-valued
solution starting from zero. Showing that this is impossible is exactly
what uniqueness means.

\medskip

The previous discussion suggests that the natural state space should be
large enough to contain signed measures arising as differences of
solutions, while still allowing one to control their behavior at
infinity. To this end, we introduce a hierarchy of spaces based on
functions with at most linear growth.

We begin with the weight function
\[
\Psi_\ell(x)
\eqdef
1+\|x\|,
\qquad x\in\R^n
\]
($\ell$ for ``linear'').

\paragraph{$\CC_\ell$: continuous functions with linear growth}

The space $\CCb(\R^n)$ of bounded continuous functions is often too
restrictive for filtering theory. For instance, the function
$x\to x$, which plays a fundamental role when considering moments of
the conditional distribution, does not belong to $\CCb(\R^n)$.

It is therefore natural to enlarge $\CCb(\R^n)$ while still retaining
some control on the behavior of functions at infinity. To this end, we
introduce the intermediate space
\[
  \CC_\ell(\R^n)
  \eqdef
  \Bigl\{
    \varphi \in \CC(\R^n) \;;\;
    \varphi/\Psi_\ell \in \CCb(\R^n)
  \Bigr\}
\]
In other words, $\CC_\ell(\R^n)$ consists of continuous functions
growing at most linearly at infinity. It is naturally endowed with the norm
\[
\|\varphi\|_{\ell}
\eqdef
\sup_{x\in\R^n}
\frac{|\varphi(x)|}{1+\|x\|}.
\]

\paragraph{$\EE_\ell$: time-dependent functions}

For a fixed horizon $T>0$, we denote by
$\EE_\ell([0,T]\times\R^n)$ the space of continuous functions
$u:[0,T]\times\R^n\to\R$ satisfying
\[
\sup_{(t,x)\in[0,T]\times\R^n}
\frac{|u(t,x)|}{1+\|x\|}
<
\infty.
\]
In other words, functions in $\EE_\ell([0,T]\times\R^n)$ are
uniformly controlled in time and exhibit at most linear growth with
respect to the spatial variable.

The space $\EE_\ell([0,T]\times\R^n)$ should be viewed primarily as a
\emph{growth space} rather than a \emph{regularity space}. Its purpose
is to control the behaviour of functions at infinity and to provide a
natural time-dependent counterpart of $\CC_\ell(\R^n)$.

It is important to note that $\EE_\ell([0,T]\times\R^n)$ is not, by
itself, the solution space of the backward equations appearing in the
uniqueness proof. Indeed, the differential operator $\LL$ involves
first- and second-order derivatives, so additional smoothness is
required. The actual solutions of the backward equation will typically
belong to
\[
\CC^{1,2}\bigl([0,T]\times\R^n\bigr)
\cap
\EE_\ell\bigl([0,T]\times\R^n\bigr),
\]
that is, they possess the regularity needed for the operator $\LL$ to
be well defined while retaining the linear-growth property encoded by
$\EE_\ell$.

Thus, $\EE_\ell([0,T]\times\R^n)$ should be regarded as the natural
ambient growth space in which sufficiently regular solutions of the
backward equation are constructed.

\paragraph{$\MM_\ell$ the state space for the weak Zakai equation}

The introduction of the larger space $\CC_\ell(\R^n)$ of test
functions has a natural counterpart on the measure side. Indeed,
functions in $\CC_\ell(\R^n)$ are no longer bounded, and therefore not
every finite signed measure can be paired with them. In order for the
duality pairing
\[
\mu(\varphi)
=
\int_{\R^n}\varphi\,\rmd\mu
\]
to be well defined for every
$\varphi\in\CC_\ell(\R^n)$, we must restrict the class of admissible
measures.

This leads naturally to the following subspace of
$\MM(\R^n)$, the space of finite signed measures on $\R^n$:
\[
\MM_\ell(\R^n)
=
\Bigl\{
\mu \in \MM(\R^n)
\;;\;
  \textstyle
  \int \Psi_\ell\,\rmd |\mu|
  =
  \int \Psi_\ell\,\rmd \mu^+
  +
  \int \Psi_\ell\,\rmd \mu^-
  < \infty
\Bigr\}.
\]

Here $|\mu|$ denotes the total variation measure of $\mu$
(\footnote{If $\mu=\mu^+-\mu^-$ is the Jordan decomposition of the signed
measure $\mu$, then $|\mu|=\mu^++\mu^-$. In particular, if $\mu$ is
positive, then $|\mu|=\mu$.}).

Equivalently, $\MM_\ell(\R^n)$ is the space of finite signed measures
with finite first moment.

The space $\MM_\ell(\R^n)$ is endowed with the weak topology induced by
$\CC_\ell(\R^n)$, i.e. a sequence
$(\mu_n)_{n\geq 1}\subset \MM_\ell(\R^n)$ converges to
$\mu\in\MM_\ell(\R^n)$ if
$\mu_n(\varphi)
\to
\mu(\varphi)$
for every test function
$\varphi\in\CC_\ell(\R^n)$.

The spaces $\CC_\ell(\R^n)$ and $\MM_\ell(\R^n)$ form a natural dual
pair. Indeed, for every
$\mu\in\MM_\ell(\R^n)$ and every
$\varphi\in\CC_\ell(\R^n)$, the pairing
$\mu(\varphi)
=
\int_{\R^n}
\varphi\,\rmd\mu
$
is finite and satisfies the estimate
$|\mu(\varphi)|
\le
\|\varphi\|_\ell\,\int \Psi_\ell\,\rmd |\mu|$.

In other words, measures in $\MM_\ell(\R^n)$ act continuously on
functions with at most linear growth. This makes
$\MM_\ell(\R^n)$ a natural state space for the uniqueness theory of the
Zakai equation. The finite first-moment condition provides precisely the
amount of integrability needed to control the action of the generator on
linearly growing test functions and to ensure that the various terms
appearing in the Zakai equation are well defined.

\paragraph{$\UU_\ell$ the class where the solution of the Zakai equation lives}

We may now finally reveal the space in which the solution of the Zakai
equation lives and in which uniqueness will be established.

Recall that the Zakai equation describes the evolution of a
measure-valued process. Since uniqueness is proved by considering the
difference of two solutions, the natural state space is
$\MM_\ell(\R^n)$ rather than the cone of positive measures.

We therefore introduce the class $\UU_\ell$ of
$\YY_t$-adapted, càdlàg, $\MM_\ell(\R^n)$-valued stochastic processes
$\mu=(\mu_t)_{t\geq 0}$
such that, for every test function
$\varphi\in\CC_\ell(\R^n)$,
\[
\tilde{\E}
\left[
\int_0^t
\bigl(\mu_s(\varphi)\bigr)^2
\,\rmd s
\right]
<
\infty,
\qquad
t>0.
\]
In other words, the process
$t\to \mu_t(\varphi)$ is required to be square-integrable on finite
time intervals for every test function of at most linear growth.

The class $\UU_\ell$ is large enough to contain the solutions arising in
the filtering problem, yet sufficiently restrictive to rule out
pathological measure-valued processes. It is within this class that the
uniqueness of the Zakai equation will be established.

\paragraph{Assumptions on the coefficients}
The spaces introduced above were designed to accommodate functions and
measures with at most linear growth. It is therefore natural to impose
regularity assumptions on the coefficients that are compatible with this
framework.

Indeed, if the coefficients were allowed to grow too rapidly at
infinity, the generator $\LL$ would no longer preserve the class of
test functions $\CC_\ell(\R^n)$, and several quantities appearing in
the Zakai equation could fail to be well defined on
$\MM_\ell(\R^n)$. Conversely, coefficients with bounded derivatives
are globally Lipschitz and grow at most linearly, which is precisely the
growth regime encoded by the spaces $\CC_\ell(\R^n)$,
$\EE_\ell([0,T]\times\R^n)$, and $\MM_\ell(\R^n)$.

For this reason, we shall work under the following regularity
assumption:

\begin{assumption}[Condition {\rm(U)}]\label{cond:U}
The components $f_i$, $a_{ij}$, and $h_i$ belong to
$\CC^2(\R^n)$. Moreover, all their partial derivatives of first and
second order belong to $\CCb(\R^n)$.
\end{assumption}

Notice that the coefficients themselves are not
assumed to be bounded. Since their first derivatives are bounded, the
functions $f$, $a$, and $h$ are globally Lipschitz and therefore
grow at most linearly.
This growth regime is precisely the one encoded by the spaces
$\CC_\ell(\R^n)$ and $\MM_\ell(\R^n)$ introduced above. In
particular, these spaces provide the natural functional framework for
the uniqueness theory developed in the sequel.

\subsection{A rigorous statement of uniqueness and its proof framework}

We can now state the main uniqueness theorem:
\begin{tcolorbox}[
 colback=gray!10,
 colframe=gray!60,
 width=\textwidth,
 boxrule=0.3pt,
 arc=2pt
]
{\vskip-0.4em\hskip-1em\color{gray!80}\sf Uniqueness of the solution of the weak Zakai equation}\\
Assume Condition~\ref{cond:U} and let
$\rho_0\in\MM_\ell(\R^n)$.
Then the weak Zakai equation \eqref{eq:zakai}
admits at most one solution in the class
$\UU_\ell$.
\end{tcolorbox}
This resuly concerns uniqueness within the class
$\UU_\ell$.
More precisely, uniqueness is understood among
$\YY_t$-adapted measure-valued processes satisfying the weak Zakai
equation together with the integrability assumptions encoded by
$\UU_\ell$.  This is the notion naturally associated with nonlinear filtering, since the filter itself is defined through conditioning with respect to the observation filtration. There are other notions of uniqueness, see Section~\ref{sec.uniqueness.notions}

The proof follows exactly the strategy outlined in
Section~\ref{sec.zakai.uniqueness.formal}. The difference is that we
can now rely on the functional framework and assumptions introduced in
the previous subsection. The purpose of the present section is
therefore to identify precisely which parts of the formal argument
require justification and to explain how the various probabilistic and
analytical ingredients fit together to yield a rigorous proof.

At a high level, three ingredients must be established rigorously:
\begin{itemize}
\item[\circnum{1}] 
      Extension of the weak Zakai equation \eqref{eq:zakai} to 
      time-dependent test functions.
\item[\circnum{2}] 
      The fact that the family of exponential martingales
      $R_t$ defined in \eqref{eq.zakai.uniqneness.R_t} is indeed a separating family of
      $\YY_T$-measurable random variables.
\item[\circnum{3}] 
      Existence and uniqueness of sufficiently regular solutions of the
      backward equation
      \eqref{eq.zakai.uniqneness.backward.pde}.
\end{itemize}
The spaces introduced previously play different roles in the argument:
\begin{itemize}
\item
$\CC_\ell$ 
provides the natural class of test
functions against which measures are evaluated.
\item
$\MM_\ell$
provides the state space for signed measures arising
as differences of solutions, and  ensures that the duality pairing
$\mu(\varphi)=\int_{\R^n}\varphi\,\rmd\mu$
is well defined for every $\varphi\in\CC_\ell$.

\item
$\EE_\ell$
controls the growth at infinity of the
time-dependent functions appearing in the backward equation.
\item
$\CC^{1,2}\cap\EE_\ell$
is the regularity class in which the backward
equation is solved.  $\CC^{1,2}$ supplies the differentiability
required by the operator $\LL$, while $\EE_\ell$ controls the
behaviour at infinity.
\item
$\UU_\ell$ specifies the class of
$\YY_t$-adapted measure-valued processes in which uniqueness is
sought.

\end{itemize}
Beyond these three main ingredients, several technical issues must also
be addressed, including the validity of It\^o's formula and the product
rule, the distinction between local martingales and true martingales,
the derivation of suitable integrability estimates, the interchange of
expectations, deterministic and stochastic integrals, and limits, as
well as the justification of the various duality pairings. Once the
estimates associated with the spaces
$\CC_\ell$, $\EE_\ell$, $\MM_\ell$, and $\UU_\ell$ are
available, these points can be handled by standard arguments from
stochastic analysis.

\medskip

Point~\circnum{1} is proved by combining a density argument and a
time-discretization argument (see \cite[Lemma 4.8]{bain2008a}). One first extends the weak Zakai equation \eqref{eq:zakai}
from $\DD(\LL)$ to $C_b^2(\R^n)$ using dominated convergence and
It\^o's isometry, and then introduces time-dependent test functions in
$C_b^{1,2}([0,T]\times\R^n)$ by freezing the time variable on a
partition of $[0,T]$ and letting the mesh size tend to zero. The
growth conditions encoded in the class $\UU_\ell$ guarantee that all
limits are justified.

\bigskip

The real difficulties are Points~\circnum{2} and~\circnum{3}.
Point~\circnum{2} is fundamentally probabilistic, since one must prove
that the family of exponential martingales forms a separating class for
$\YY_T$-measurable random variables. Point~\circnum{3} is
fundamentally analytical, as it requires existence, uniqueness, and
regularity results for the backward equation
\eqref{eq.zakai.uniqneness.backward.pde}. 

This latter task is made more subtle
by two features of the filtering model. First, the coefficients are not
assumed to be bounded, which motivates the introduction of the growth
spaces $\CC_\ell$, $\EE_\ell$, and $\MM_\ell$. 

Second, the diffusion operator need not be uniformly elliptic,
preventing a direct application of the classical existence,
uniqueness, and regularity results for uniformly parabolic equations.
This places the analysis within the more delicate framework of
degenerate parabolic PDEs and explains the appeal to more advanced
results from parabolic PDE theory (see \cite{friedman1964a}).

The backward equation belongs to the class of parabolic partial
differential equations. More precisely, it is an evolution equation
involving a first-order derivative in time together with the
second-order differential operator associated with the diffusion
process. Such equations arise naturally as the PDE counterparts of
stochastic differential equations.

\subsubsection{The separating family of exponential martingales  $R_t$ }

At first sight, the introduction of the complex exponential martingale
\eqref{eq.zakai.uniqneness.R_t} may appear somewhat mysterious. Its
origin is in fact closely related to a classical idea from Fourier
analysis.

Indeed, the complex exponentials
$x \to e^{i\lambda x}$, with $\lambda\in\R$,
form a separating family of real functions. Roughly speaking, if a real function
$f$ satisfies
$\int_{\R} f(x)\,e^{i\lambda x}\,\rmd x = 0$,
for all $\lambda\in\R$,
then $f=0$. This uniqueness property lies at the heart of Fourier
analysis and of the theory of characteristic functions in probability.

The same idea is used here in an infinite-dimensional setting. The
processes
\[
t\longmapsto
\exp\left(
i\int_0^t r_s^\top\,\rmd Y_s
\right),
\]
where $r$ ranges over a sufficiently rich class of deterministic
functions, play the role of Fourier modes on the observation path space.
After the appropriate It\^o correction, i.e. adding $\frac12\int_0^t \|r_s\|^2\,\rmd s$ inside the exponential, which does not affect the separation argument, they become the exponential
martingales $(R_t)_{t\ge0}$.

A first important feature is that $R_t$ is {\bfseries bounded}. Indeed, since
$|e^{iz}|=1$,
the oscillatory factor
$\exp\left(
i\int_0^t r_s^\top\,\rmd Y_s
\right)
$
has modulus one. Consequently,
\[
|R_t|
=
\exp\left(
\frac12\int_0^t \|r_s\|^2\,\rmd s
\right).
\]
If $r$ is bounded, then
\[
|R_t|
\le
\exp\left(
\frac12\,T\,\|r\|_\infty^2
\right),
\qquad 0\le t\le T,
\]
so that $R_t$ is uniformly bounded and therefore integrable.

A second important feature is that $R_t$ is a {\bfseries martingale}. This allows
one to apply It\^o's formula, integration by parts, and conditional
expectation arguments in a rigorous way. The martingale property is
precisely what makes $R_t$ a convenient testing device in the
uniqueness proof.

Most importantly, the family of random variables $R_T$ is rich enough
to {\bfseries separate $\YY_T$-measu\-rable random variables}, just as the
family $(e^{i\lambda x})_{\lambda\in\R}$ separates real functions. Once one has established
\eqref{eq.zakai.uniqneness.separation.cond}
for every admissible choice of $r$, this separating property implies
that $\nu_T(\psi)=0$, and therefore that the two candidate solutions
must coincide.

In this sense, the exponential martingales $R_t$ may be viewed as the
analogue of Fourier exponentials on Wiener space. Their role is not to
describe the filtering problem itself, but rather to provide a
sufficiently rich family of test random variables with which uniqueness
can be established.

\subsubsection{The backward dual equation}

The third ingredient is the analysis of the backward equation \eqref{eq.zakai.uniqneness.backward.pde}.
This equation appears naturally in the duality argument where its role is to
cancel the drift term arising in the dynamics \eqref{eq.zakai.uniqneness.R_t.nu_t} of the quantity $R_t\,\nu_t(\phi_t)$.

The terminal condition $\varphi$ is left arbitrary. By varying
$\varphi$ over a sufficiently rich class of test functions, one
obtains a sufficiently rich family of backward solutions capable of
probing the measure $\nu_T$ in all directions.

The main question is whether this backward equation admits a unique
solution possessing enough regularity for the duality argument to be
valid.

The desired solution class is $\CC^{1,2}\cap\EE_\ell$.
The factor $\CC^{1,2}$ is required because the operator $\LL$
contains first-order derivatives in time and second-order derivatives
in space. The factor $\EE_\ell$ ensures that the solution exhibits
the growth properties compatible with the spaces
$\CC_\ell$ and~$\MM_\ell$.

Two sources of difficulty arise simultaneously.

First, the coefficients of the signal--observation system are not
assumed to be bounded. This is important, since excluding unbounded
coefficients would rule out many classical models, including the
linear-Gaussian case underlying the Kalman--Bucy filter. One therefore
cannot work exclusively with bounded test functions and finite measures,
and suitable growth conditions must be introduced instead.

Second, the diffusion operator need not be uniformly elliptic. As a
result, the backward equation does not automatically fall within the
simplest framework of uniformly parabolic PDEs, and one cannot directly
appeal to the strongest existence, uniqueness, and regularity results
from classical parabolic PDE theory.

These two features are largely responsible for the technical complexity
of the uniqueness theory.

\paragraph{Unbounded coefficients}

Throughout these notes, the coefficients of the state--observation
system are assumed to satisfy regularity conditions ensuring existence
and uniqueness of the underlying stochastic differential equations.
However, they are not assumed to be bounded.

More precisely, boundedness of the first derivatives typically implies
at most linear growth of the coefficients. Consequently, quantities
such as $f(x)$, $a(x)$, $h(x)$
may become arbitrarily large as $\|x\|\to\infty$.

This feature has several important consequences.

First, one cannot restrict attention to bounded test functions. The
natural backward equations generate solutions whose growth must be
controlled rather than eliminated. This motivates the introduction of
the growth spaces
\[
\CC_\ell(\R^n)
\qquad\text{and}\qquad
\EE_\ell([0,T]\times\R^n),
\]
which allow at most linear growth in the spatial variable.

Second, the measures appearing in the Zakai equation must possess
sufficient moments to make the duality pairing
\[
\mu(\varphi)
=
\int_{\R^n}\varphi(x)\,\mu(\rmd x)
\]
well defined for all admissible test functions
$\varphi\in\CC_\ell(\R^n)$. This leads naturally to the space $\MM_\ell(\R^n)$.

Finally, many integrability estimates appearing throughout the proof
must be formulated in terms of moment bounds rather than uniform
bounds. In this sense, the spaces
$\CC_\ell$, $\EE_\ell$, and $\MM_\ell$ should be viewed as the
natural framework in which the linear growth of the coefficients can be
handled.

\paragraph{Lack of uniform ellipticity.}

The second difficulty concerns the backward equation appearing in the
duality argument. Recall that the operator
\[
\LL
=
\sum_{i=1}^{n}
f_i(x)\,\partial_i
+
\frac12
\sum_{i,j=1}^{n}
a_{ij}(x)\,\partial_{ij}
\]
is the infinitesimal generator of the signal process. The associated
backward equation
\[
\partial_t \phi
+
\LL\phi
+
i\,h^\top r_t\,\phi
=
0
\]
belongs to the class of parabolic partial differential equations.

The analysis of such equations is particularly well understood when the
diffusion matrix is \emph{uniformly elliptic}, that is, when there
exists a constant $\lambda>0$ such that
\[
\xi^\top a(x)\,\xi
\ge
\lambda \|\xi\|^2,
\qquad
x,\xi\in\R^n.
\]
From a probabilistic viewpoint, this condition means that randomness is
injected in every direction of the state space. The resulting diffusion
acts as a smoothing mechanism, spreading mass and regularizing
solutions. Consequently, the backward equation falls within the
classical theory of uniformly parabolic PDEs, where existence,
uniqueness, and regularity follow from a powerful collection of tools,
including maximum principles, Schauder estimates, fundamental
solutions, and semigroup methods. In particular, one can often obtain
solutions in
\[
\CC^{1,2}\bigl([0,T]\times\R^n\bigr)
\]
under relatively mild assumptions on the coefficients.

However, Condition~\ref{cond:U} does not impose uniform ellipticity.
It only requires smoothness of the coefficients together with bounded
first and second derivatives. The diffusion matrix may therefore be
degenerate, and some directions of the state space may receive little
or no noise. As a consequence, the backward equation does not fall
directly within the simplest framework of uniformly parabolic
equations, and the classical theory cannot be applied directly.

The analysis is nevertheless still possible. Bain and Crisan appeal to
results from Friedman's monograph
\emph{Partial Differential Equations of Parabolic Type}
\cite{friedman1964a}, which provides existence, uniqueness, and
regularity results under assumptions substantially weaker than uniform
ellipticity. In particular, uniqueness is obtained through a
maximum-principle argument (Lemma~4.12), while existence follows from
deeper results on linear parabolic equations. These results provide the
PDE ingredient required by the duality method and ultimately make the
uniqueness proof rigorous.

\subsection{Further remarks on uniqueness}

\subsection{Alternative proof}
\label{sec.uniqueness.alternative}

The duality method presented above is not the only available approach
to uniqueness. In fact, when the signal noise and the observation noise
are independent, considerably simpler proofs can be obtained.

This is precisely the framework considered throughout these notes.
Nevertheless, we have chosen to present the duality argument because it
extends naturally to more general settings, including correlated
observation noise. The resulting proof therefore captures the ideas
underlying the general theory rather than the simplifications specific
to the independent-noise case.

The basic observation is that, under independence of the noises, the
weak Zakai equation can be rewritten in a mild form involving the
Markov semigroup $(P_t)_{t\ge0}$ associated with the signal process:
\[
\rho_t(\varphi)
=
\rho_0(P_t\varphi)
+
\int_0^t
\rho_s\!\left(
P_{t-s}(\varphi h^\top)
\right)
\,\rmd Y_s .
\]

This formula may be viewed as the stochastic analogue of the
variation-of-constants formula for ordinary differential equations.
Instead of introducing a backward adjoint equation, one propagates the
test function directly through the semigroup $P_t$.

If $\rho^1$ and $\rho^2$ are two solutions of the Zakai equation
with the same initial condition and
\[
\nu_t=\rho_t^1-\rho_t^2,
\]
then
\[
\nu_t(\varphi)
=
\int_0^t
\nu_s\!\left(
P_{t-s}(\varphi h^\top)
\right)
\,\rmd Y_s .
\]
The problem is thus reduced to a stochastic integral equation.
Combining It\^o's isometry with standard moment estimates and a
Grönwall argument, one can show that
$\nu_t(\varphi)=0$ for every admissible test function $\varphi$,
and hence that $\nu_t=0$.

Historically, uniqueness proofs of this type were developed by several
authors, notably Jacques Szpirglas, who exploited semigroup techniques
and mild formulations of the Zakai equation. Such arguments are often
technically simpler than the Bensoussan duality method in the
independent-noise setting.

The price paid for this simplification is a loss of generality. The
mild representation above relies heavily on the independence of the
signal and observation noises and does not extend naturally to the
correlated case. By contrast, the duality method adapts much more
readily to different observation models and has therefore become one of
the standard approaches in the general theory of nonlinear filtering.

For this reason, although simpler proofs are available in the setting
considered here, we have chosen to present the duality argument.
Beyond its role in proving uniqueness, it illustrates one of the most
important ideas in stochastic analysis: the interaction between forward
stochastic evolutions and backward adjoint equations.

\subsection{Which notion of uniqueness?}
\label{sec.uniqueness.notions}

Before leaving the uniqueness question, it is worth clarifying the
precise meaning of the uniqueness result proved above. Unlike ordinary
differential equations, stochastic equations may admit several distinct
notions of uniqueness, depending on the class of admissible solutions.

The uniqueness theorem established in these notes should be understood
as follows:

\begin{center}
\emph{There exists at most one $\YY_t$-adapted solution belonging to
the class $\UU_\ell$.}
\end{center}

This is the natural notion of uniqueness in nonlinear filtering.
Indeed, the filter is intended to represent the observer's knowledge of
the hidden state based solely on the information contained in the
observation process. Since the filtering problem itself is formulated
through conditional expectations with respect to the observation
filtration
\[
\YY_t
=
\sigma(Y_s,\;0\le s\le t),
\]
both the normalized filter $\pi_t$ and the unnormalized filter
$\rho_t$ are automatically $\YY_t$-adapted. Adaptedness is
therefore not merely a technical assumption; it is an intrinsic part of
the filtering problem.

From a broader perspective, one may distinguish between several levels
of uniqueness. At the filtering level, one seeks uniqueness among
adapted solutions that genuinely describe the evolution of information
available to the observer. This is the viewpoint adopted throughout
these notes. A stronger question consists in viewing the Zakai equation
as a stochastic partial differential equation in its own right and
asking whether it admits a unique solution among all measure-valued
stochastic processes satisfying the equation, regardless of any
adaptedness requirement. Such results belong to the theory of
measure-valued SPDEs and typically require different techniques and
additional assumptions.

The distinction is analogous to that encountered in the theory of
stochastic differential equations, where one may study uniqueness among
adapted solutions on a given probability space or uniqueness in law.
The uniqueness proof presented here belongs to the former category: it
establishes uniqueness within the class of $\YY_t$-adapted
measure-valued solutions satisfying the growth conditions encoded in
$\UU_\ell$.

For the purposes of nonlinear filtering, this is precisely the notion
that matters. It guarantees that the Zakai equation uniquely
characterizes the evolution of the unnormalized conditional
distribution associated with the observation process.

\subsection{Historical and bibliographical comments}
\label{rem.uniqueness.history}

The uniqueness proof presented in these notes follows closely the
approach developed in \textcite{bain2008a}, which remains one of the
standard modern references on nonlinear filtering. Our objective has
not been to improve upon the original proof, but rather to make
explicit some of the ideas, motivations, and intermediate steps that
are often left implicit in more condensed treatments.

Historically, the proof combines ingredients originating from two
distinct mathematical traditions. The probabilistic component relies on
martingale methods and stochastic exponentials. In particular, the
family
\[
R_t
=
\exp\!\left(
i\int_0^t r_s^\top\,\rmd Y_s
+
\frac12
\int_0^t \|r_s\|^2\,\rmd s
\right)
\]
may be viewed as providing Fourier modes on Wiener space, allowing one
to separate $\YY_t$-measurable random variables in much the same way
that ordinary Fourier exponentials separate functions on Euclidean
spaces.

The analytical component is rooted in the classical duality principle
of partial differential equations. The central idea is to pair a
forward evolution equation with a solution of the corresponding
adjoint backward equation in such a way that the drift terms cancel.
This philosophy appears throughout the theory of Markov processes,
Kolmogorov equations, Fokker--Planck equations, and semigroup theory,
and long predates nonlinear filtering itself.

In the filtering literature, the resulting uniqueness argument is often
referred to as the \emph{Bensoussan method}. What is specific to
nonlinear filtering is the combination of the classical duality
argument with stochastic exponentials and observation filtrations. The
proof presented above is therefore best viewed as an interaction
between two complementary ideas: a probabilistic Fourier analysis on
path space and a PDE duality argument.

The treatment of the backward equation relies heavily on results from
the theory of parabolic partial differential equations. For this
reason, \textcite{bain2008a} repeatedly appeals to Friedman's
monograph
\emph{Partial Differential Equations of Parabolic Type}
\cite{friedman1964a}. Its appearance in the proof is a useful reminder
that nonlinear filtering is not solely a probabilistic theory. Many of
its central results ultimately depend on deep analytical properties of
parabolic PDEs.

Finally, as discussed in
Section~\ref{sec.uniqueness.alternative}, alternative uniqueness proofs
are available under stronger assumptions. In particular, when the
signal and observation noises are independent, semigroup methods and
mild formulations often lead to considerably shorter arguments.
Nevertheless, the duality method remains one of the most flexible
approaches and extends naturally to more general observation models.

\section{The Kushner--Stratonovich equation}

{\it
\color{black!45!SecColor!60!white}
The Zakai equation provides a linear stochastic evolution equation for the
\emph{unnormalized} conditional distribution of the signal. While this
representation is mathematically convenient, the quantity of primary
interest in filtering is the \emph{normalized} conditional distribution,
which describes the posterior law of the state given the observations.

Historically, the Kushner--Stratonovich equation was derived before the
Zakai equation and was originally introduced as the fundamental evolution
equation for the nonlinear filter. However, because it is nonlinear, its
mathematical analysis is considerably more involved. The later discovery
of the Zakai equation provided a linear description of the filtering
problem in terms of an unnormalized measure, which is often easier to
study both theoretically and numerically. For this reason, it has become
standard to introduce the Zakai equation first and then recover the
Kushner--Stratonovich equation through normalization.

The purpose of this section is to derive the evolution equation satisfied
by the normalized filter. Starting from the Zakai equation and using the
Kallianpur--Striebel normalization formula, we obtain the
Kushner--Stratonovich equation, which is the fundamental nonlinear
stochastic equation of continuous-time filtering. We then introduce the
innovation process, discuss the interpretation of the resulting dynamics,
and conclude with a formal derivation of its strong (density) form.

}

\bigskip

We now derive the Kushner--Stratonovich equation from the Zakai equation.

The Zakai equation governs the \emph{unnormalized} filter $(\rho_t)_{t\ge0}$,
whereas the nonlinear filter of interest is the \emph{normalized}
conditional distribution $(\pi_t)_{t\ge0}$. These two objects are related
by the Kallianpur--Striebel formula:
\[
\pi_t(\varphi)=\frac{\rho_t(\varphi)}{\rho_t(1)}.
\]

Thus, the Kushner--Stratonovich equation is obtained by normalizing the
Zakai equation.

\subsection{Derivation via Itô's formula}

Starting from the Zakai equation \eqref{eq:zakai}, we first apply it to
the constant test function $\varphi \equiv 1$, which yields
\[
\rho_t(1)
=
\rho_0(1)
+
\int_0^t \rho_s(h^\top)\,\rmd Y_s.
\]
We then consider the normalized filter defined by
\[
\pi_t(\varphi)
=
\frac{\rho_t(\varphi)}{\rho_t(1)}.
\]
To derive its dynamics, we apply Itô's formula to this quotient. Using
the stochastic differential equations satisfied by $\rho_t(\varphi)$ and
$\rho_t(1)$, together with the quadratic covariation
\[
\rmd \langle \rho(\varphi), \rho(1)\rangle_t
=
\rho_t(\varphi h^\top)\,\rho_t(h)\,\rmd t,
\]
a straightforward computation leads to
\begin{tcolorbox}[
 colback=gray!10,colframe=gray!60,width=\textwidth,boxrule=0.3pt,arc=2pt
]
{\vskip-0.4em\hskip-1em\color{gray!80}\sf the Kushner--Stratonovich equation}
\begin{align}
\label{eq:ks}
\pi_t(\varphi)
=
\pi_0(\varphi)
+
\int_0^t \pi_s(\LL\varphi)\,\rmd s
+
\int_0^t
  \bigl(\pi_s(\varphi h^\top)-\pi_s(\varphi)\;\pi_s(h^\top)\bigr)
  \;
  \bigl(\rmd Y_s-\pi_s(h)\,\rmd s\bigr).
\end{align}
for all $t\geq 0$ and $\varphi \in \Cc^2(\R^n)$.
\end{tcolorbox}
In contrast with the Zakai equation, the Kushner--Stratonovich equation is nonlinear due to the normalization. 

As expected, in the absence of observations ($h \equiv 0$), \eqref{eq:ks} reduces to the weak form of the Fokker--Planck equation \eqref{eq:fp}.

\subsection{Innovation process}

We now introduce
\begin{tcolorbox}[
 colback=gray!10,colframe=gray!60,width=\textwidth,boxrule=0.3pt,arc=2pt
]
{\vskip-0.4em\hskip-1em\color{gray!80}\sf the innovation process}
\[
\rmd \nu_t \eqdef \rmd Y_t-\pi_t(h)\,\rmd t\,.
\]
\end{tcolorbox}
It represents the new information
contained in the observations after subtracting the predicted drift.

In terms of the innovation process, the equation takes the form
\[
\pi_t(\varphi)
=
\pi_0(\varphi)
+
\int_0^t \pi_s(\LL\varphi)\,\rmd s
+
\int_0^t
\Big(
\pi_s(\varphi h^\top)-\pi_s(\varphi)\pi_s(h^\top)
\Big)\,\rmd \nu_s.
\]

\medskip

The term
$\pi_s(\varphi h^\top)-\pi_s(\varphi)\pi_s(h^\top)$
is the conditional covariance under $\pi_s$:
\[
\mathrm{Cov}_{\pi_s}(\varphi, h)
=
\pi_s(\varphi h^\top)-\pi_s(\varphi)\pi_s(h^\top).
\]

With this notation, the equation becomes
\[
\pi_t(\varphi)
=
\pi_0(\varphi)
+
\int_0^t \pi_s(\LL\varphi)\,\rmd s
+
\int_0^t
\mathrm{Cov}_{\pi_s}(\varphi, h)\,\rmd \nu_s.
\]

\subsection{Interpretation}

The Kushner--Stratonovich equation has a natural interpretation in terms
of a continuous-time prediction--correction mechanism.

The term
\[
\int_0^t \pi_s(\LL\varphi)\,\rmd s
\]
corresponds to the \emph{prediction step}, driven by the intrinsic
dynamics of the state process.

The term
\[
\int_0^t \mathrm{Cov}_{\pi_s}(\varphi, h)\,\rmd \nu_s
\]
corresponds to the \emph{correction step}, incorporating the new
information brought by the observations.

The innovation process
\[
\rmd \nu_t = \rmd Y_t - \pi_t(h)\,\rmd t
\]
represents the unexpected part of the observation, that is, the difference
between what is observed and what was predicted.

\medskip

The factor $\mathrm{Cov}_{\pi_s}(\varphi, h)$ determines how strongly the
observation affects the estimate: it measures the dependence between the
quantity of interest $\varphi(X_s)$ and the observation function
$h(X_s)$ under the current conditional distribution.

\begin{itemize}
\item If $\varphi(X_s)$ and $h(X_s)$ are weakly correlated, the correction is small.
\item If they are strongly correlated, the observation has a significant impact.
\end{itemize}

\medskip

Thus, the Zakai and Kushner--Stratonovich equations provide two
complementary descriptions of the filtering problem: the former is linear
but unnormalized, while the latter is normalized but nonlinear.

\subparagraph{Bayesian interpretation}
{\it 
The Kallianpur--Striebel formula can be viewed as a continuous-time
version of Bayes' rule.

The unnormalized filter $\rho_t$ corresponds to the product of prior and
likelihood, whose evolution is linear (Zakai equation). The normalized filter $\pi_t$ is obtained from the same dynamics, completed by a renormalization, and represents the posterior distribution (Kushner--Stratonovich equation).
 This normalization is precisely what introduces the
nonlinearity in the filtering dynamics.

Thus, the Zakai equation describes the linear evolution of
``likelihood $\times$ prior'', while the Kushner--Stratonovich equation
describes the evolution of the posterior distribution.
}

\subsection{Strong form of the Kushner--Stratonovich equation}

We now give a \emph{formal} derivation of the strong form of the
Kushner--Stratonovich equation.

Recall that the normalized filter satisfies
\[
\pi_t(\varphi)
=
\frac{\rho_t(\varphi)}{\rho_t(1)}.
\]
Assume that, for each $t\ge0$, $\pi_t$ admits a density with respect to
Lebesgue measure, say
\[
\pi_t(\rmd x)=p_t(x)\,\rmd x.
\]
Then, for every test function $\varphi$,
\[
\pi_t(\varphi)
=
\int_{\R^n} \varphi(x)\,p_t(x)\,\rmd x
=
\crochet{\pi_t,\varphi}
=
(p_t,\varphi).
\]

Using the Kushner--Stratonovich equation in weak form, we have
\begin{align*}
(p_t,\varphi)
=
(p_0,\varphi)
+
\int_0^t (p_s,\LL\varphi)\,\rmd s
+
\int_0^t
\Big(
(p_s,\varphi h^\top)
-
(p_s,\varphi)\,(p_s,h^\top)
\Big)
\,\rmd \nu_s,
\end{align*}
where $\rmd \nu_s = \rmd Y_s - \pi_s(h)\,\rmd s$ is the innovation process.

As in the Zakai case, we introduce the adjoint operator $\LL^*$ defined by
\[
(u,\LL\varphi)=(\LL^*u,\varphi).
\]
Using this relation, we obtain
\[
(p_s,\LL\varphi)
=
(\LL^*p_s,\varphi).
\]

Moreover,
\[
(p_s,\varphi h^\top)
=
(p_s\,h^\top,\varphi),
\qquad
(p_s,\varphi)\,(p_s,h^\top)
=
(p_s\,\pi_s(h)^\top,\varphi).
\]
Therefore,
\begin{align*}
(p_t,\varphi)
&=
(p_0,\varphi)
+
\int_0^t (\LL^*p_s,\varphi)\,\rmd s
+
\int_0^t
\bigl(
p_s\,(h-\pi_s(h))^\top,\varphi
\bigr)
\,\rmd \nu_s.
\end{align*}

Since this identity holds for all test functions $\varphi$, we are led,
formally, to
\begin{tcolorbox}[
 colback=gray!10,colframe=gray!60,width=\textwidth,boxrule=0.3pt,arc=2pt
]
{\vskip-0.4em\hskip-1em\color{gray!80}\sf the (strong) Kushner--Stratonovich equation}
\begin{align}
\label{eq:ks.strong}
\rmd p_t(x)
&=
\LL^* p_t(x)\,\rmd t
+
p_t(x)\,\bigl(h(x)-\pi_t(h)\bigr)^\top\,\rmd \nu_t,
\end{align}
for $t>0$ and $x\in\R^n$, with initial condition $p_0$ equal to
the density of the law of $X_0$.
\end{tcolorbox}

Equation~\eqref{eq:ks.strong} is a stochastic partial differential equation (SPDE) satisfied by the density of the normalized filter.

\medskip

This derivation is purely formal. A rigorous treatment requires additional
regularity assumptions ensuring the existence of a density, the validity
of the integration by parts argument, and a suitable interpretation of
\eqref{eq:ks.strong} (typically in a weak or variational sense).

\medskip

As expected, in the absence of observations ($h \equiv 0$), \eqref{eq:ks.strong} reduces to the strong form of the Fokker--Planck equation \eqref{eq:fp.strong}.

\section{Finite-dimensional filters}

{\it
\color{black!45!SecColor!60!white}

As we have seen, nonlinear filtering is intrinsically
infinite-dimensional. A natural question therefore arises: are there
situations in which the nonlinear filter $\pi_t$ can be represented
within a known family of probability distributions parameterized by a
finite number of variables? In such a case, the filtering problem would
reduce to the estimation of the finite-dimensional parameters
characterizing the conditional distribution, rather than the evolution
of the distribution itself. This question, known as the finite-dimensional filtering problem, is precisely the topic of this section.
}

\bigskip

The Zakai and Kushner--Stratonovich equations show that nonlinear
filtering is naturally an infinite-dimensional problem. Indeed, the
object of interest is the conditional distribution
\[
\pi_t
=
\law(X_t\mid\YY_t),
\]
whose evolution is governed by a measure-valued stochastic equation.
Except in very special situations, the full conditional distribution
must therefore be propagated in time.

A natural question is whether this infinite-dimensional dynamics can
sometimes be reduced to a finite-dimensional system. More precisely, one
may ask whether the conditional distribution can be represented by a
finite number of parameters whose evolution is governed by a closed
system of differential equations.

The most celebrated example is provided by the Kalman--Bucy filter. In
the linear-Gaussian setting, the conditional distribution remains
Gaussian for all times and is therefore completely characterized by its
mean and covariance matrix. The infinite-dimensional filtering problem
then collapses to a finite-dimensional dynamical system.

\subsection{The linear-Gaussian model and the Kalman--Bucy filter}

Consider the state-observation system
\begin{align}
  \rmd X_t
  &=
  A \,X_t\,\rmd t + B\,\rmd W_t\,,
  \qquad
  X_0
\sim
\NN(\bar X_0,Q_0)\,,
\\
  \rmd Y_t
  &=
  C \,X_t\,\rmd t
  +
  D\,\rmd V_t\,,
\end{align}
where $W_t$ and $V_t$ are independent Brownian motions and also independent from the initial condition $X_0$.

Because linear transformations preserve Gaussianity and because Gaussian
distributions are stable under  conditioning, the conditional
distribution remains Gaussian:
\[
\pi_t
=
\law(X_t|\YY_t)
=
\NN(m_t,P_t).
\]
and hence depends only on the two first conditional moments of $\pi_t$:
\begin{align*}
  m_t
  &\eqdef
  \E[X_t|\YY_t]\,,
  &
  P_t
  &\eqdef
  \mathrm{Cov}[X_t|\YY_t]
\end{align*}
Consequently, instead of propagating the full conditional distribution,
it suffices to propagate the conditional mean and covariance matrix, leading to the Kalman-Bucy filter\wiki{https://en.wikipedia.org/wiki/Kalman_filter\#Kalman–Bucy_filter}:
\begin{align}
\label{eq.kalman.hat.X}
  \rmd m_t
  &=
  A \,m_t\,\rmd t 
  + K_t\; \bigl( \rmd Y_t - C \,m_t\,\rmd t \bigr)\,,
  & m_0=\bar X_{0}\,,
\\
\label{eq.kalman.P}
  \dot  P_t
  &=
  A P_t
  + P_t A^\top
  + BB^\top - P_t C^\top (DD^\top)^{-1} C P_t\,,
  & P_{0} = Q_{0}\,,
\end{align}
where $K_{t}$ is the  Kalman gain defined as:
\[
  K_t \eqdef P_t \,C^\top \,(DD^\top)^{-1}\,.
\]
A remarkable feature of the Kalman filter is that the covariance matrix
$P_t$ evolves according to a deterministic Riccati equation, whereas
the conditional mean $m_t$ satisfies a stochastic equation driven by
the observations. This reflects the linear-Gaussian structure of the
model: although $m_t$ depends on the particular observation path, the
conditional covariance is entirely determined by the system parameters
and the initial covariance. In particular, $P_t$ can be computed
\emph{offline}, before any observations are available, by solving the
Riccati equation once and for all. The observations affect the estimate
itself, but not the evolution of the estimation uncertainty.

The innovation process is given by
\[
\nu_t
=
Y_t
-
\int_0^t C \,m_s\,\rmd s.
\]
As in the general nonlinear filtering theory, $(\nu_t)_{t\ge0}$ is a
Brownian motion with respect to the observation filtration. It measures
the discrepancy between the incoming observations and the prediction
provided by the current estimate.

\medskip

The key ingredient behind the Kalman-Bucy filter is a classical result on
Gaussian vectors: if
\[
  \left(\begin{smallmatrix}
    X \\ Y
  \end{smallmatrix}\right)
  \sim
  \NN
  \left(
  \left(\begin{smallmatrix}
  \bar X\\
  \bar Y
  \end{smallmatrix}\right)
  \;,\;
  \left(\begin{smallmatrix}
  \Sigma_{XX} & \Sigma_{XY}\\
  \Sigma_{YX} & \Sigma_{YY}
  \end{smallmatrix}\right)
  \right).
\]
is jointly Gaussian, then the conditional distribution of $X$ given
$Y$ is again Gaussian, conditional mean $\hat X$ and covariance $P$ given by:
\begin{align*}
  \hat X
  &=
  \bar X
  +
  \Sigma_{XY}\Sigma_{YY}^{-1}\,\bigl(Y-\bar Y\bigr)
  &
  P
  &=
  \Sigma_{XX}
  -
  \Sigma_{XY}\Sigma_{YY}^{-1}\Sigma_{YX}
\end{align*}
and $K
\eqdef
\Sigma_{XY}\Sigma_{YY}^{-1}$,
which is precisely the structure appearing in the correction step of
the Kalman filter.

\subsection{A finite-dimensional realization of the filter}

The Kalman filter is called a \emph{finite-dimensional filter} because
the infinite-dimensional object $\pi_t$ is completely characterized by
the finite-dimensional pair
\[
(m_t,P_t).
\]

Equivalently, the filtering dynamics preserve the finite-dimensional
manifold of Gaussian distributions
\[
\Bigl\{
\NN(m,P)
:
m\in\R^n,
\;
P\in\mathbb S_n^+
\Bigr\}\,,
\]
where $\mathbb S_n^+$ denotes the set of $n\times n$ symmetric positive semidefinite matrices.

The Kalman filter can therefore be viewed as an exact finite-dimensional
realization of the Kushner--Stratonovich equation.

Indeed, the conditional mean and covariance can be expressed directly in
terms of the filter:
\[
m_t
=
\pi_t(x\mapsto x)
=
\E[X_t\mid\YY_t],
\]
and
\[
P_t
=
\pi_t\Bigl(
x\mapsto (x-m_t)(x-m_t)^\top
\Bigr)
=
\E\left[
(X_t-m_t)(X_t-m_t)^\top
\mid
\YY_t
\right].
\]

In the Gaussian case, these two moments completely determine the
conditional distribution. Consequently, the infinite-dimensional process
$(\pi_t)_{t\ge0}$ may be replaced by the finite-dimensional process
$(m_t,P_t)_{t\ge0}$.

\subsection{Other finite-dimensional filters}
\label{sec.finite.dimensional.filters}

A natural question is whether other finite-dimensional filters exist.
More precisely, one may ask whether there are nonlinear
state-observation systems for which the conditional distribution remains
confined to a finite-dimensional family of probability measures.

This question motivated a substantial body of research initiated by
Vic Bene\v{s}, Roger Brockett, Sanjoy Mitter, Michiel Hazewinkel, Fred Daum and others. Several
remarkable examples have been discovered, most notably the Bene\v{s}
filter, in which a specific nonlinear structure preserves a
finite-dimensional exponential family of distributions.

Nevertheless, such examples are rare and highly structured. For generic
nonlinear systems, finite-dimensional realizations do not exist and the
filtering problem remains intrinsically infinite-dimensional. In
particular, a remarkable generic non-existence result was established by
Mireille Chaleyat-Maurel and Dominique Michel. Roughly speaking,
their work shows that finite-dimensional filters are exceptional objects
and cannot be expected to arise for generic nonlinear diffusion models.

This result reinforces the central role of the Zakai and
Kushner--Stratonovich equations. Rather than seeking a finite collection
of sufficient statistics, these equations describe directly the
evolution of the full conditional distribution, viewed as an
infinite-dimensional dynamical object.

From this perspective, the Kalman filter is not a different theory but a
remarkable special case of nonlinear filtering. The infinite-dimensional
dynamics of the conditional distribution happen to preserve the
finite-dimensional manifold of Gaussian measures, allowing the filter to
be represented exactly by a finite number of parameters. In general,
however, no such reduction is available, and the measure-valued
equations of filtering theory remain unavoidable.

\section{Outro: a review of assumptions}

The filtering problem is formulated in a weak (measure-valued) sense: the objects of interest (e.g. $\pi_t$, $\rho_t$) act on test functions rather than being described pointwise. Because the dynamics involve the infinitesimal generator $\LL$, we require function spaces that both characterize measures and are sufficiently regular for Itô's formula and the action of $\LL$. We now  review the different classes of test functions used:
\begin{description}
\item[\normalfont\color{SecColor}Bounded measurable functions.]
$\pi_t$ is defined by
$\pi_t(\varphi)=\E[\varphi(X_t)| \YY_t]$
for all $\varphi \in \Bb(\R^n)$.
This class characterizes measures but is too large for differential
operations.
\item[\normalfont\color{SecColor}Continuous bounded functions.]
The space $\Cb(\R^n)$ is measure-determining\footnote{A class $\CC$ is measure-determining if $\mu(\varphi)=\nu(\varphi)$ for all $\varphi\in\CC$ implies $\mu=\nu$.}
and provides a convenient topological setting, but lacks sufficient
regularity.
\item[\normalfont\color{SecColor}Twice continuously differentiable functions.]
We work with $\Cc^2(\R^n)$ (or functions with sufficient decay), since:
\fenumi\ 
Itô's formula applies to $\varphi \in \Cc^2(\R^n)$.
\fenumii\ 
$\LL$ is naturally defined on $\Cc^2(\R^n)$.
\fenumiii\ 
Decay ensures integration by parts without boundary terms, yielding
$\LL^*$ and the strong (SPDE) formulation.
\fenumiv\ 
$\Cc^2(\R^n)$ is measure-determining.
\fenumv\ 
Although $\Cc^2(\R^n)$ is not the full domain $\DD(\LL)$ of $\LL$, it forms a convenient core on which the generator is well defined and all required computations can be carried out.
\end{description}

We also briefly review the main {\color{SecColor}assumptions on the coefficients} beyond measurability:
\begin{description}
\item[\normalfont\color{SecColor}State equation.]
We assume that $f$ and $g$ are globally Lipschitz. This ensures a unique
strong solution and implies at most linear growth, yielding standard
moment bounds. If $X_0 \in L^2$ (Assumption~\ref{ass:model-finite-moments}),
then $X_t \in L^2$ for all $t>0$, which suffices for defining conditional
expectations. The derivation of the filtering equations requires an additional
integrability assumption, namely that $X_0$ admits a finite third moment
(Assumption~\ref{ass:model-finite-moments.bis}).
\item[\normalfont\color{SecColor}Observation function.]
We assume that $h$ has at most linear growth, ensuring that the
observation process and the likelihood integrals are well defined.
\item[\normalfont\color{SecColor}Novikov condition and integrability.]
To ensure that the likelihood process $(Z_t)_{t\geq 0}$ is a martingale, we use
Novikov's condition \eqref{eq.novikov}. It holds automatically if $h$ is
bounded, but this excludes the linear/Gaussian case.
If $h$ has linear growth, Novikov's condition requires stronger
integrability of $(X_t)_{t\geq 0}$. Since $|h(X_s)|^2 \lesssim 1+|X_s|^2$, it follows
from exponential moment estimates. These are ensured if $X_0$ admits an
exponential moment (Assumption~\ref{ass:model-exp-moments}), which covers
the Gaussian setting.
\end{description}

\section{Scope and references}

In these notes, we restrict ourselves to a simple diffusion framework in
order to highlight the main ideas of nonlinear filtering. Several
extensions are not considered.

First, we assume that the coefficients are time-homogeneous, i.e.\ of the
form $f(x)$, $g(x)$, and $h(x)$. Extensions to time-dependent
coefficients, such as $f(t,x)$, can be treated in a similar manner, at
the expense of slightly heavier notation.

We also assume that the coefficients do not depend on the observation
process. Models with such a dependence can be treated, but they lead to a
slightly different formulation of the filtering problem.

Throughout, we suppose that $(W_t)_{t\ge0}$ and $(V_t)_{t\ge0}$ are
independent. This simplifies the structure of the filtering problem.
More general models with correlated noises can also be studied, but they
do not change the main ideas and complicate the exposition.

We further restrict ourselves to the case where the observation noise has
constant covariance. Allowing the covariance of $(V_t)_{t\ge0}$ to depend
on the state leads to significantly more involved formulations and raises
nontrivial issues in the analysis of the filtering equations; see  \cite{crisan2011a}.

Finally, we do not consider models beyond the diffusion setting, such as
piecewise deterministic Markov processes (PDMP), Lévy-driven systems, or
diffusions with jumps. 

\bigskip

The literature on stochastic filtering is vast. We briefly comment on the main references used in these notes.

Several general books on stochastic calculus provide the necessary background, in particular those of 
\cite{karatzas1991b}, \cite{oksendal2013a}, and \cite{revuz2005a}.

The main monographs dedicated to nonlinear filtering for diffusion processes are those of \cite{bain2008a}, and \cite{xiong2008a}, which give a clear and systematic treatment of the subject.
Among the older references, the collective volume by \cite{mitter1983a} and
the Saint-Flour lecture notes by \cite{pardoux1991a} remain important
sources on nonlinear filtering.

A comprehensive reference is the handbook edited by \cite{crisan2011a}, which gathers contributions from several authors and covers many aspects of nonlinear filtering.

The book by Alain \cite{bensoussan1992b} occupies a special place: it develops the theory under assumptions that are sometimes different from the standard diffusion framework, and contains tools that allow one to treat cases including the linear Gaussian setting (Kalman filtering).

A recent contribution by \cite{crisan2026} studies the uniqueness of the filtering equations in a general framework, in particular allowing for observation noise with state-dependent covariance.

The inclusion of \cite{friedman1964a} serves as a reminder that the
theory of nonlinear filtering is not purely probabilistic. In
particular, the uniqueness proof presented in these notes ultimately
relies on classical results from the theory of parabolic partial
differential equations.

We have chosen not to cite the original articles in which nonlinear
filtering for diffusion processes was developed; the interested reader
may find them through the references listed above. Moreover, virtually
all of the material presented in these notes can be found, in one form
or another, in the excellent monograph of Bain and Crisan
\cite{bain2008a}.

\clearpage
\nocite{*}
\addcontentsline{toc}{section}{References}
\printbibliography

\appendix

\end{document}